\documentclass[english,ruled]{article}
\usepackage[T1]{fontenc}
\usepackage[latin9]{inputenc}
\usepackage{verbatim}
\usepackage[algo2e,ruled]{algorithm2e}
\usepackage{amsmath}
\usepackage{amsthm}
\usepackage{amssymb}
\usepackage{graphicx}
\usepackage{xcolor}
\usepackage{enumitem}
\makeatletter
\usepackage[toc,page,header]{appendix}
\usepackage{minitoc}
\usepackage[accepted]{icml2023}
\usepackage{comment}
\usepackage{ifthen}
\newboolean{doublecolumn}
\setboolean{doublecolumn}{true}
\usepackage{multirow}
\usepackage{pgfplots}
\usepackage{pdflscape}

\newcommand{\condsep}{\ifthenelse{\boolean{doublecolumn}}{}{\\}}

\renewcommand \thepart{}
\renewcommand \partname{}

\theoremstyle{plain}
\newtheorem{lem}{\protect\lemmaname}
\theoremstyle{remark}
\newtheorem{rem}{\protect\remarkname}
\theoremstyle{plain}
\newtheorem{thm}{\protect\theoremname}
\theoremstyle{plain}

\providecommand{\corollaryname}{Corollary}
\theoremstyle{plain}

\usepackage{babel}
\providecommand{\propositionname}{Proposition}
\usepackage{hyperref}       
\usepackage{url}            
\usepackage{booktabs}       
\usepackage{amsfonts}       
\usepackage{nicefrac}       
\usepackage{authblk}
\usepackage{optidef}
\usepackage{array}

\makeatother

\usepackage{babel}

\usepackage[resetlabels]{multibib}

\newcites{app}{References in the Appendix}

\providecommand{\lemmaname}{Lemma}
\providecommand{\remarkname}{Remark}
\providecommand{\theoremname}{Theorem}

\icmltitlerunning{Fast online algorithms for offline LPs}

\begin{document}

\global\long\def\inprod#1#2{\left\langle #1,#2\right\rangle }%

\global\long\def\inner#1#2{\langle#1,#2\rangle}%

\global\long\def\binner#1#2{\big\langle#1,#2\big\rangle}%

\global\long\def\Binner#1#2{\Big\langle#1,#2\Big\rangle}%

\global\long\def\norm#1{\lVert#1\rVert}%

\global\long\def\bnorm#1{\big\Vert#1\big\Vert}%

\global\long\def\Bnorm#1{\Big\Vert#1\Big\Vert}%

\global\long\def\setnorm#1{\Vert#1\Vert_{-}}%

\global\long\def\bsetnorm#1{\big\Vert#1\big\Vert_{-}}%

\global\long\def\Bsetnorm#1{\Big\Vert#1\Big\Vert_{-}}%


\global\long\def\brbra#1{\big(#1\big)}%
\global\long\def\Brbra#1{\Big(#1\Big)}%
\global\long\def\rbra#1{(#1)}%


\global\long\def\sbra#1{[#1]}%
\global\long\def\bsbra#1{\big[#1\big]}%
\global\long\def\Bsbra#1{\Big[#1\Big]}%


\global\long\def\cbra#1{\{#1\}}%
\global\long\def\bcbra#1{\big\{#1\big\}}%
\global\long\def\Bcbra#1{\Big\{#1\Big\}}%
\global\long\def\vertiii#1{\left\vert \kern-0.25ex  \left\vert \kern-0.25ex  \left\vert #1\right\vert \kern-0.25ex  \right\vert \kern-0.25ex  \right\vert }%

\global\long\def\matr#1{\bm{#1}}%

\global\long\def\til#1{\tilde{#1}}%

\global\long\def\wtil#1{\widetilde{#1}}%

\global\long\def\wh#1{\widehat{#1}}%

\global\long\def\mcal#1{\mathcal{#1}}%

\global\long\def\mbb#1{\mathbb{#1}}%

\global\long\def\mtt#1{\mathtt{#1}}%

\global\long\def\ttt#1{\texttt{#1}}%

\global\long\def\dtxt{\textrm{d}}%

\global\long\def\bignorm#1{\bigl\Vert#1\bigr\Vert}%

\global\long\def\Bignorm#1{\Bigl\Vert#1\Bigr\Vert}%

\global\long\def\rmn#1#2{\mathbb{R}^{#1\times#2}}%

\global\long\def\deri#1#2{\frac{d#1}{d#2}}%

\global\long\def\pderi#1#2{\frac{\partial#1}{\partial#2}}%

\global\long\def\limk{\lim_{k\rightarrow\infty}}%

\global\long\def\trans{\textrm{T}}%

\global\long\def\onebf{\mathbf{1}}%

\global\long\def\zerobf{\mathbf{0}}%

\global\long\def\zero{\bm{0}}%


\global\long\def\Euc{\mathrm{E}}%

\global\long\def\Expe{\mathbb{E}}%

\global\long\def\rank{\mathrm{rank}}%

\global\long\def\range{\mathrm{range}}%

\global\long\def\diam{\mathrm{diam}}%

\global\long\def\epi{\mathrm{epi} }%

\global\long\def\relint{\mathrm{relint} }%

\global\long\def\inte{\operatornamewithlimits{int}}%

\global\long\def\cov{\mathrm{Cov}}%

\global\long\def\argmin{\operatornamewithlimits{arg\,min}}%

\global\long\def\argmax{\operatornamewithlimits{arg\,max}}%

\global\long\def\tr{\operatornamewithlimits{tr}}%

\global\long\def\dis{\operatornamewithlimits{dist}}%

\global\long\def\sign{\operatornamewithlimits{sign}}%

\global\long\def\prob{\mathrm{Prob}}%

\global\long\def\spans{\textrm{span}}%

\global\long\def\st{\operatornamewithlimits{s.t.}}%

\global\long\def\dom{\mathrm{dom}}%

\global\long\def\prox{\mathrm{prox}}%

\global\long\def\for{\mathrm{for}}%

\global\long\def\diag{\mathrm{diag}}%

\global\long\def\and{\mathrm{and}}%

\global\long\def\st{\mathrm{s.t.}}%

\global\long\def\dist{\mathrm{dist}}%

\global\long\def\Var{\operatornamewithlimits{Var}}%

\global\long\def\raw{\rightarrow}%

\global\long\def\law{\leftarrow}%

\global\long\def\Raw{\Rightarrow}%

\global\long\def\Law{\Leftarrow}%

\global\long\def\vep{\varepsilon}%

\global\long\def\dom{\operatornamewithlimits{dom}}%

\global\long\def\tsum{{\textstyle {\sum}}}%

\global\long\def\Cbb{\mathbb{C}}%

\global\long\def\Ebb{\mathbb{E}}%

\global\long\def\Fbb{\mathbb{F}}%

\global\long\def\Nbb{\mathbb{N}}%

\global\long\def\Rbb{\mathbb{R}}%

\global\long\def\extR{\widebar{\mathbb{R}}}%

\global\long\def\Pbb{\mathbb{P}}%

\global\long\def\Mrm{\mathrm{M}}%

\global\long\def\Acal{\mathcal{A}}%

\global\long\def\Bcal{\mathcal{B}}%

\global\long\def\Ccal{\mathcal{C}}%

\global\long\def\Dcal{\mathcal{D}}%

\global\long\def\Ecal{\mathcal{E}}%

\global\long\def\Fcal{\mathcal{F}}%

\global\long\def\Gcal{\mathcal{G}}%

\global\long\def\Hcal{\mathcal{H}}%

\global\long\def\Ical{\mathcal{I}}%

\global\long\def\Kcal{\mathcal{K}}%

\global\long\def\Lcal{\mathcal{L}}%

\global\long\def\Mcal{\mathcal{M}}%

\global\long\def\Ncal{\mathcal{N}}%

\global\long\def\Ocal{\mathcal{O}}%

\global\long\def\Pcal{\mathcal{P}}%

\global\long\def\Scal{\mathcal{S}}%

\global\long\def\Tcal{\mathcal{T}}%

\global\long\def\Wcal{\mathcal{W}}%

\global\long\def\Xcal{\mathcal{X}}%

\global\long\def\Ycal{\mathcal{Y}}%

\global\long\def\Zcal{\mathcal{Z}}%

\global\long\def\i{i}%


\global\long\def\abf{\mathbf{a}}%

\global\long\def\bbf{\mathbf{b}}%

\global\long\def\cbf{\mathbf{c}}%

\global\long\def\dbf{\mathbf{d}}%

\global\long\def\fbf{\mathbf{f}}%

\global\long\def\lambf{\bm{\lambda}}%

\global\long\def\alphabf{\bm{\alpha}}%

\global\long\def\sigmabf{\bm{\sigma}}%

\global\long\def\thetabf{\bm{\theta}}%

\global\long\def\deltabf{\bm{\delta}}%

\global\long\def\sbf{\mathbf{s}}%

\global\long\def\lbf{\mathbf{l}}%

\global\long\def\ubf{\mathbf{u}}%

\global\long\def\vbf{\mathbf{v}}%

\global\long\def\wbf{\mathbf{w}}%

\global\long\def\xbf{\mathbf{x}}%

\global\long\def\ybf{\mathbf{y}}%

\global\long\def\zbf{\mathbf{z}}%

\global\long\def\Abf{\mathbf{A}}%

\global\long\def\Ubf{\mathbf{U}}%

\global\long\def\Pbf{\mathbf{P}}%

\global\long\def\Ibf{\mathbf{I}}%

\global\long\def\Ebf{\mathbf{E}}%

\global\long\def\Mbf{\mathbf{M}}%

\global\long\def\Qbf{\mathbf{Q}}%

\global\long\def\Lbf{\mathbf{L}}%

\global\long\def\Pbf{\mathbf{P}}%


\global\long\def\abm{\bm{a}}%

\global\long\def\bbm{\bm{b}}%

\global\long\def\cbm{\bm{c}}%

\global\long\def\dbm{\bm{d}}%

\global\long\def\ebm{\bm{e}}%

\global\long\def\fbm{\bm{f}}%

\global\long\def\gbm{\bm{g}}%

\global\long\def\hbm{\bm{h}}%

\global\long\def\pbm{\bm{p}}%

\global\long\def\qbm{\bm{q}}%

\global\long\def\rbm{\bm{r}}%

\global\long\def\sbm{\bm{s}}%

\global\long\def\tbm{\bm{t}}%

\global\long\def\ubm{\bm{u}}%

\global\long\def\vbm{\bm{v}}%

\global\long\def\wbm{\bm{w}}%

\global\long\def\xbm{\bm{x}}%

\global\long\def\ybm{\bm{y}}%

\global\long\def\zbm{\bm{z}}%

\global\long\def\Abm{\bm{A}}%

\global\long\def\Bbm{\bm{B}}%

\global\long\def\Cbm{\bm{C}}%

\global\long\def\Dbm{\bm{D}}%

\global\long\def\Ebm{\bm{E}}%

\global\long\def\Fbm{\bm{F}}%

\global\long\def\Gbm{\bm{G}}%

\global\long\def\Hbm{\bm{H}}%

\global\long\def\Ibm{\bm{I}}%

\global\long\def\Jbm{\bm{J}}%

\global\long\def\Lbm{\bm{L}}%

\global\long\def\Obm{\bm{O}}%

\global\long\def\Pbm{\bm{P}}%

\global\long\def\Qbm{\bm{Q}}%

\global\long\def\Rbm{\bm{R}}%

\global\long\def\Ubm{\bm{U}}%

\global\long\def\Vbm{\bm{V}}%

\global\long\def\Wbm{\bm{W}}%

\global\long\def\Xbm{\bm{X}}%

\global\long\def\Ybm{\bm{Y}}%

\global\long\def\Zbm{\bm{Z}}%

\global\long\def\lambm{\bm{\lambda}}%

\global\long\def\alphabm{\bm{\alpha}}%

\global\long\def\albm{\bm{\alpha}}%

\global\long\def\taubm{\bm{\tau}}%

\global\long\def\mubm{\bm{\mu}}%

\global\long\def\yrm{\mathrm{y}}%

\global\long\def\assign{\coloneqq}%
\global\long\def\tmstrong#1{{\bf #1}}%
\global\long\def\tmem#1{#1}%
\global\long\def\tmop#1{#1}%

\global\long\def\tmverbatim#1{#1}%
\global\long\def\revised#1{{#1}}%
\global\long\def\changed#1{\textcolor{blue}{#1}}%

\newtheorem{definition}{Definition}
\newcommand{\TODO}[1]{\textbf{\textcolor{red}{TODO: #1}}}

\twocolumn[
\icmltitle{Solving Linear Programs with Fast Online Learning Algorithms}



\icmlsetsymbol{equal}{*}

\begin{icmlauthorlist}
\icmlauthor{Wenzhi Gao}{icme}
\icmlauthor{Dongdong Ge}{sjtu}
\icmlauthor{Chunlin Sun}{icme}
\icmlauthor{Yinyu Ye}{mse}
\end{icmlauthorlist}
\icmlaffiliation{icme}{ICME, Stanford University}
\icmlaffiliation{mse}{Management Science and Engineering Stanford University}
\icmlaffiliation{sjtu}{Antai College of Economics and Management, Shanghai Jiaotong University}

\icmlcorrespondingauthor{Wenzhi Gao}{gwz@stanford.edu}

\icmlkeywords{Online linear programming, Large-scale linear programming, Column generation}
\vskip 0.3in
]



\printAffiliationsAndNotice{} 

\newcommand{\tmtextbf}[1]{\text{{\bfseries{#1}}}}
\newcommand{\tma}{\mathbf{a}}
\newcommand{\tmb}{\mathbf{b}}
\newcommand{\tmc}{\mathbf{c}}
\newcommand{\tmd}{\mathbf{d}}
\newcommand{\e}{\mathbf{e}}
\newcommand{\1}{\textbf{1}}
\newcommand{\A}{\mathbf{A}}
\newcommand{\x}{\mathbf{x}}
\newcommand{\y}{\mathbf{y}}
\newcommand{\s}{\mathbf{s}}
\newcommand{\tmu}{\mathbf{u}}
\newcommand{\0}{\textbf{0}}
\newcommand{\dl}{\underline{d}}
\newcommand{\du}{\bar{d}}
\newcommand{\w}{\mathbf{w}}
\newcommand{\myfrac}[2]{\tfrac{#1}{#2}}
\newcommand{\condsmall}{\ifthenelse{\boolean{doublecolumn}}{\small}{}}

\setlength{\parindent}{0pt}

\begin{abstract}
  This paper presents fast first-order methods for solving linear programs (LPs) approximately. We adapt online linear programming algorithms to offline LPs and obtain algorithms that avoid any matrix multiplication. We also introduce a variable-duplication technique that copies each variable $K$ times and reduces the optimality gap and constraint violation by a factor of $\sqrt{K}$. Furthermore, we show how online algorithms can be effectively integrated into sifting, a column generation scheme for large-scale LPs. Numerical experiments demonstrate that our methods can serve as either an approximate direct solver, or an initialization subroutine for exact LP solving.
\end{abstract}
\section{Introduction}
First-order methods for large-scale linear programs are receiving increasing
attention as problem sizes grow beyond the capacity of traditional simplex
{\cite{dantzig2016linear}} and interior point methods {\cite{ye2011interior}}.
When a low or moderate-accuracy solution is required,
first-order methods often outperform simplex and interior point solvers
as they only involve matrix-vector multiplication and are free of matrix
decomposition.
A recent practice {\cite{applegate2021practical}}
demonstrates that fine-tuned first-order methods are also capable of solving
LPs of millions of variables to high accuracy.
However, as {\cite{nesterov2012efficiency}} suggests,
when the problem size further reaches trillions of variables or more, even storing the problem data in memory becomes prohibitive, not to mention
performing matrix-vector multiplications. 
As a result, most first-order methods either resort to distributed architecture
{\cite{basu2020eclipse, fountoulakis2016performance}} or rely on sparsity {\cite{nesterov2012efficiency}}.\condsep

In contrast to the LPs discussed above, which we call ``offline'' LPs, there is another type of LPs called ``online'' LPs
{\cite{agrawal2014dynamic, 10.48550/arxiv.1909.05499}}.
They come from revenue management {\cite{talluri2004theory}} and 
model the situation of making sequential decisions on the customer orders that are described by LP columns.
Here an LP column denotes a column in the constraint matrix together with its corresponding objective value.
Due to online nature of the problem, generally an algorithm for online LPs
\tmtextbf{1)}. only accesses one LP column at a time; \tmtextbf{2)}. makes
decisions instantly. Through making decisions, the behavior of online algorithms gets refined and improved,
which is why they are also known as ``online learning''.
\condsep

In the line of online LP literature, state-of-the-art online LP
algorithms build on offline LPs as a subroutine
{\cite{10.48550/arxiv.2210.07996, 10.48550/arxiv.1909.05499, 10.48550/arxiv.2209.00399}},
and it seems that \textit{online LP unilaterally \textbf{benefits from} offline LP}. It is therefore
natural to ask:
\begin{center}
{\textit{Does online LP benefit offline LP?}} 
\end{center}
In this paper, we answer this question positively by exploiting a recent research direction on fast online algorithms 
{\cite{10.48550/arxiv.2003.02513, 10.48550/arxiv.2011.10124}}. 
Unlike the offline LP-based online algorithms, these methods build on stochastic first-order oracles, make decisions and refinements in
linear time, and simultaneously estimate primal and dual
solutions in a single pass through the problem data. Also, they achieve provable optimality guarantees under minimal assumptions on the problem data. All these
motivate us to study if these time-efficient methods help offline LP solving. In particular, we
investigate two popular online algorithms: explicit subgradient {\cite{10.48550/arxiv.2003.02513}} and implicit proximal point
{\cite{kulis2010implicit}}.

\paragraph{First-order Methods for Linear Programming}

First-order methods for LPs are recently developed and presented as
optimization software. Some successful attempts include \texttt{ABIP} \cite{lin2021admm,deng2022new}, \texttt{COSMO}
\cite{garstka2021cosmo}, \texttt{ECLIPSE} \cite{basu2020eclipse}, \texttt{PDLP} \cite{applegate2021practical} and \texttt{SCS}
\cite{o2016conic}. The listed algorithms target large-scale LPs, and
the major computational bottleneck lies in matrix-vector multiplication (sometimes one-time matrix decomposition).
Compared to the aforementioned works, our work is more in line with (block) subgradient-based or coordinate descent-based approaches
{\cite{nesterov2012efficiency, nesterov2014subgradient}}. What differentiates our method is that they {\cite{nesterov2012efficiency, nesterov2014subgradient}} have to update primal and dual variables sequentially, while our proposed method simultaneously updates both solutions.

\paragraph{Simple and Fast Online LP}

A simple and fast online LP algorithm is recently proposed and analyzed by
{\cite{10.48550/arxiv.2003.02513}}. It exploits the finite-sum structure of the LP dual problem and 
  runs stochastic subgradient method in the dual space, while simultaneously estimating the primal solution by optimality condition. 
  Under mild assumptions, given an LP of $m$ constraints and $n$ variables, the authors show that
their algorithm achieves ${\mathcal{O}} \left( \sqrt{n} \log n + m \sqrt{n} \right)$ constraint violation and an ${\mathcal{O}} (m \sqrt{n} )$ optimality gap (also known as regret in online learning literature) under the random permutation setting that we exploit in this paper. 
Concurrent works {\cite{balseiro2020dual, 10.48550/arxiv.2011.10124}}
further generalize the results to nonlinear objective using a mirror descent
framework and achieve $\mathcal{O} (m +  \sqrt{m n} )$ gap with no constraint violation under a weaker setting.
We also note that \cite{10.48550/arxiv.2011.10124} covers an adversarial setting that is more general than random permutation. However, the optimality gap is too pessimistic to be applied for offline LPs, and therefore we exploit the setting of {\cite{10.48550/arxiv.2003.02513}}.

\paragraph{Contributions}

In this paper, we show that online LP algorithms can benefit offline LP both theoretically and practically. Specifically, 
\begin{itemize}[leftmargin=8pt]
  \item We investigate how two online learning algorithms perform when they are used to solve offline LPs.
  These two methods, known as explicit subgradient and implicit proximal point in the online learning literature, are
  free of any matrix-vector multiplication through sequential online access to LP columns.
  We provide theoretical analyses of the optimality gap and
  constraint violation leveraging the analysis from online LP literature.
  For both algorithms, we obtain  an ${\mathcal{O}} \left(
   m \log n + \sqrt{n} \log n + \sqrt{m n} \right)$ optimality gap and ${\mathcal{O}} \left(m+ \sqrt{m n} \right)$ constraint violation, which improves the results from {\cite{10.48550/arxiv.2003.02513}}. Moreover, to enhance the practical performance of the online algorithms and to make them usable for offline LPs, we take advantage of the offline LP setting and propose a variable duplication
  scheme which further reduces the gap and violation by a factor of $\sqrt{K}$, where $K$ is the number of copies of each variable.
  \item For exact LP solving, we identify that our algorithms can perfectly fit into sifting \cite{bixby1992very}, an LP column generation framework targeting problems with far more columns than constraints. We show online algorithms provide both an
  initial basis guess and an approximate dual solution. Numerical experiments
  show that our online algorithms can accelerate sifting by
  providing a good initialization and by stabilizing the dual solution.
\end{itemize}

\begin{table*}[h]
  \caption{Previous optimality gap and constraint violation bounds. Ex: explicit update; Im: implicit proximal point update; Duplicate: \textbf{Algorithm \ref{alg:online-dup}}; $\rho$: optimality gap; $v$: constraint violation. We note that the Duplicate setting is only available for offline LPs.
}
  \centering
  \small
  \begin{tabular}{cccc}
    \toprule
    Work & Algorithm & Stochastic input ($v$ \& $\rho$)  & Random permutation ($v$ \& $\rho$)\\
    \midrule
     \cite{10.48550/arxiv.2003.02513} & Ex. Subgrad & $\mathcal{O} \left( m + m \sqrt{n}
    \right)$ \& $\mathcal{O} \left( m \sqrt{n} \right)$  & $\mathcal{O}
    \left( m \sqrt{n} \right)$ \& $\mathcal{O} \left( \sqrt{n} \log n + m
    \sqrt{n} \right)$\\
    \midrule
    \cite{10.48550/arxiv.2011.10124} & Ex. Mirror & $\mathcal{O} \left( m + \sqrt{m n} \right)$
    \& $0$ & ---\\
    \midrule
    & Ex. Subgrad & $\mathcal{O} \left( m + \sqrt{m n} \right)$
    \& $\mathcal{O} \left( \sqrt{m n} \right)$ & $\mathcal{O} \left( m +
    \sqrt{m n} \right)$ \& $\mathcal{O} ( m \log n + \sqrt{n} \log n + \sqrt{m n}
    )$\\
    This work & Im. Proximal  & $\mathcal{O} \left( m + \sqrt{m n} \right)$ \&
    $\mathcal{O} \left( \sqrt{m n} \right)$ & $\mathcal{O} \left( m +
    \sqrt{m n} \right)$ \& $\mathcal{O} \left( m \log n + \sqrt{n} \log n + \sqrt{m n}
    \right)$\\
    & Ex/Im. Duplicate & $\mathcal{O} ( \frac{m}{K} + \sqrt\frac{{m
    n}}{{K}} )$ \& $\mathcal{O} ( \sqrt\frac{{m
    n}}{{K}} )$ & $\mathcal{O} ( \frac{m}{K} + \sqrt\frac{{m
    n}}{{K}} )$ \& $\mathcal{O} ( \myfrac{m \log n}{K}+ \sqrt{\myfrac{n}{K}} \log n + \sqrt{\myfrac{m n}{K}} )$\\
    \bottomrule
  \end{tabular}
\end{table*}

\paragraph{Other Related Work}

Disregarding the complexity of solving subproblems, offline LP-based algorithms
achieve state-of-the-art theoretical guarantees for online LPs. For example,
{\cite{10.48550/arxiv.2210.07996, 10.48550/arxiv.2209.00399, bray2019logarithmic}} all achieve $\mathcal{O}{(\log n)}$ optimality gap (up to a $\log \log $ term) under different LP non-degeneracy assumptions.\condsep

Going beyond the online LP problem, which is classified as online convex optimization with
constraint (\texttt{OCOwC)}, there is substantial research on general
online convex optimization {\cite{hazan2016introduction}}. While most online algorithms are subgradient-based, there is also another line of research based on the implicit proximal point update {\cite{kulis2010implicit, chen2022implicit}}. 
In this paper, we choose the implicit update as one of the fast algorithms.

\paragraph{Outline of the Paper} This paper is organized as follows. In \textbf{Section \ref{sec:setup}}, we introduce the basic problem setup and assumptions; In \textbf{Section \ref{sec:theory}}, we provide theoretical analyses of the performance of our algorithms; \textbf{Section \ref{sec:sifting}} briefly discusses the sifting procedure for large-scale LPs and how online algorithms can be potentially combined with this framework. Last we conduct numerical experiments in \textbf{Section \ref{sec:experiments}} to verify our theory and to show the practical efficiency of our proposed methods. In the appendix, we further elaborate on some practical aspects of our methods.
\section{Problem Setup and Assumptions} \label{sec:setup}

\paragraph{Notations}

Throughout the paper, we use $\| \cdot \|$ to denote the Euclidean norm and
$\langle \cdot, \cdot \rangle$ to denote the Euclidean inner product. Bold
letters $\A, \tma$ are used to denote matrices and vectors; the
subdifferential of a convex function $f$ is denoted by $\partial f (x)$ and $f' (x)
\in \partial f (x)$ is called a subgradient; we use $[\cdot]_+ \assign \max
\{ \cdot, 0 \}$ to denote the element-wise positive part function; $\mathbb{I} \{
\cdot \}$ denotes the 0-1 indicator function and $\delta_S$ is the indicator function of set $S$. 
Unless specified, we use $\mathbb{E} [\cdot]$ to
denote expectation taken over a certain permutation of LP columns. We use $\tma_{[i
: j]}$ to denote indexing of vector $\tma$ from the $i$-th to $j$-th
coordinate, or indexing of a matrix from the $i$-th to $j$-th column.

\subsection{Linear Programming and Duality}
Given a linear program of $m$ constraints and $n$
variables
\begin{align}
  \max_{\x \in \mathcal{F}_p}\quad & \langle \tmc, \x \rangle, \nonumber
\end{align}
its dual problem is given by
\begin{align}
  	  \min_{(\y, \s) \in \mathcal{F}_d} \quad &  \langle \tmb, \y \rangle +
\langle \1, \s \rangle \nonumber
\end{align}
where the primal and dual feasible sets are denoted by
\begin{align}
	\mathcal{F}_p := &\left\{ \x : \A \x \leq \tmb, \0 \leq \x \leq \1 \right\}\nonumber\\
  \mathcal{F}_d := &\left\{ ( \y, \s
  ) \geq \0 : \A^{\top} \y + \s \geq \tmc, \s \geq \0 \right\},\nonumber
\end{align}
with $\A \in \mathbb{R}^{m \times n}, \tmb \in \mathbb{R}^m, \tmc \in
\mathbb{R}^n$ and $\1 = (1, \ldots, 1)^{\top}$ denotes the column vector of
all ones. By the same argument as in
{\cite{10.48550/arxiv.2003.02513}}, we remove the dual variable $\s$ and
rewrite the dual problem with a finite-sum formulation.
\begin{eqnarray}
  \min_{\y \geq \0}~ \myfrac{1}{n} \small\sum_{j = 1}^n  \langle \tmd, \y
  \rangle + [ c_j - \langle \tma_j, \y \rangle ]_+
 =: \myfrac{1}{n} \sum_{j = 1}^n f( \y, j)
\end{eqnarray}
where $\tmd = \tmb / n$ and $\tma_j$ denotes the $j$-th column of $\A$.
Here, with a slight abuse of notation, we denote $f(\y, j)=\tmd, \y
  \rangle + [ c_j - \langle \tma_j, \y \rangle ]_+$ as the stochastic function associated with the LP column $(c_j,\tma_j)$ for all $j$ without the explicit mention of $c$ and $\tma$ in the function arguments.
Given the optimal dual solution $\y^*$, the optimality conditions tell us
\begin{align} \label{eqn:kkt}
x_j^{\ast}
\in \left\{ \begin{array}{ll}
  \{ 0 \}, & c_j < \langle \tma_j, \y^{\ast} \rangle\\
  {}[0, 1], & c_j = \langle \tma_j, \y^{\ast} \rangle\\
  \{ 1 \}, & c_j > \langle \tma_j, \y^{\ast} \rangle
\end{array} \right.	
\end{align}
Taking into account this finite-sum structure of the dual problem, together with the close
relevance between primal and dual solutions, it is appealing to apply first-order
stochastic algorithms to the dual problem and simultaneously estimate primal solutions using the relation \eqref{eqn:kkt}. This is indeed what ``simple and fast online algorithms'' does.

\subsection{Simple and Fast Online Algorithms}
Now, we are ready to introduce online algorithms of our interest in \textbf{Algorithm \ref{alg:online}}. 

\begin{algorithm}
\caption{Fast online algorithms for LP \label{alg:online}}	
\KwIn{$\y^0, \{f(\y, j)\}$ from LP data $\A, \tmb, \tmc$}
\For{k = \rm{1 to $n$ }}{
	  \hspace{-3pt}{Choose $i_k \in [n]$}
      \begin{eqnarray}
        \hspace{-8pt}\y^{k + 1} = \argmin_{\y \geq \0} ~ \{ f_{\y^k} ( \y, k
        ) + \myfrac{1}{2\gamma_k} \| \y - \y^k \|^2 \} \label{eqn:dual-update}
      \end{eqnarray}
      \hspace{-8pt}Estimate $\hat{x}^{i_k}$ based on \eqref{eqn:dual-update}.
}
\KwOut{$\hat{\x}$}
\end{algorithm}
In each iteration, let $\y^k$ be the current estimation of the dual solution. We choose $i_k\in[n]$, approxiamte $f(\y, i_k)$ by some function $f_{\y^k}(\y, i_k)$, which will be defined in \textbf{Section \ref{sec:theory}}, and perform stochastic proximal updates on the dual variable $\y$.
At the end of each iteration, we use the current dual information, combined with \eqref{eqn:kkt} to estimate the primal solution $\hat{x}^{i_k}$.
Since in each iteration, only one LP column participates in the update, the cost of each iteration is very low, thereby making the methods  ``simple and fast''.

\subsection{Assumptions}
We make the following assumptions across the paper.
\begin{enumerate}[leftmargin=30pt, label=\textbf{A\arabic*:},ref=\rm{\textbf{A\arabic*}}]
  \item $\max_{i} d_i = \bar{d} \geq \min_i d_i = \dl > 0 \quad(\tmd = \tmb / n)$ \label{A1}
  \item $\| \tma_j \|_{\infty} \leq \bar{a}$ and $| c_j | \leq
  \bar{c}$~ for any $j \in [n]$ \label{A2}
  \item For any dual solution $\y \geq \0$, no more than $m$ columns
  satisfy $c_j = \langle \tma_j, \y \rangle$ \label{A3}
\end{enumerate}
\begin{rem} The above assumptions are standard in the literature of online LP \cite{10.48550/arxiv.2003.02513, 10.48550/arxiv.2011.10124, 10.48550/arxiv.1909.05499}, and can also be easily satisfied by most offline LPs.
  {\ref{A1}} is satisfied by a wide range of LPs such as multi-knapsack,
  set-covering, and online resource allocation. Even if $b_i = 0$, we can perturb it to make the assumption hold; Bounds in \ref{A2} are only used for
  analysis and $\bar{a}, \bar{c}$ can be computed by a single pass through data.
  \ref{A3} could be satisfied by an arbitrarily small perturbation of $\tmb$ \cite{megiddo1989varepsilon, agrawal2014dynamic}. 
\end{rem}

\subsection{Performance Measure}
We use optimality gap and constraint
violation to measure the quality of a given primal solution $\hat{\x}$.
\begin{align}
	\rho ( \hat{\x} ) &\assign \max_{\x \in \mathcal{F}_p}~~\langle
   \tmc, \x \rangle - \langle \tmc, \hat{\x} \rangle\\
   v ( \hat{\x} ) &\assign \| [ \A \hat{\x} - \tmb ]_+ \|
\end{align}

With all the tools in hand, we are ready to present several realizations of \textbf{Algorithm \ref{alg:online}} and to analyze their performance for offline LPs.

\section{Fast Online Algorithms for Offline LPs} \label{sec:theory}

In this section, we provide two realizations of \textbf{Algorithm \ref{alg:online}} and analyze their expected optimality gap and constraint violation.  We start by specifying the choice of $i_k$ in \textbf{Algorithm \ref{alg:online}}. Although stochastic input, namely sampling $i_k$ from $[n]$ randomly with replacement is a feasible option, a problem is that, with this method, some LP columns $(c_j,\tma_j)$ might not be traversed within $\mathcal{O}(n)$ iterations, and therefore we cannot get a solution estimate for those columns. Therefore, random permutation, or choosing $i_k$ by sampling without replacement from $[n]$, is a better choice for offline LP's context. In other words, we can safely stop after $n$ iterations and ensure that all the columns are associated with their estimated primal values.

To avoid the heavy notations from permutation, without loss of generality we assume instead of sampling $i_k$ from $[n]$, we permute the columns of the offline LP, so that we can simply let $i_k = k$ and get the same theoretical results. 

\begin{rem}
	As is shown in \cite{10.48550/arxiv.2003.02513}, once we prove the result for the random permutation setting, the analysis can be directly applied to the stochastic input setting with a slightly better bound for the optimality gap.
\end{rem}

\subsection{Online Explicit Update}
\label{sec:theory-explicit}

In this section, we analyze the performance of the online explicit subgradient update in the
offline setting. The (Sub)gadient-based update approximates $f ( \y, k )$ by a linear function
\begin{align}
  f_{\y^k} ( \y, k ) & = \langle f' ( \y^k, k ),
  \y - \y^k \rangle \label{eqn:explicit} \\
  & =\langle \tmd - \tma_k \mathbb{I} \{ c_k > \langle
  \tma_k, \y^k \rangle \}, \y - \y^k \rangle, \nonumber
\end{align}
where we take $\tmd - \tma_k \mathbb{I} \{ c_k > \langle \tma_k, \y^k
\rangle \} \in \partial f ( \y^k, k )$ and estimate the
primal solution by $x^k =\mathbb{I} \{ c_k > \langle \tma_k, \y^k
\rangle \}$. Specifically, the dual update \eqref{eqn:dual-update} is given in the closed form
$$\y^{k+1} = [\y^{k} - \gamma (\tmd - \tma_k x^k)]_+. $$
With this dual updating rule, our \textbf{Algorithm \ref{alg:online}} becomes the method in {\cite{10.48550/arxiv.2003.02513}}. Compared with {\cite{10.48550/arxiv.2003.02513}}, we provide a sharper analysis to achieve a better trade-off between the optimality gap and constraint violation.

\begin{lem} \label{lem:explicit-bound}
  {Under assumptions {\ref{A1}} to \ref{A3}, if we take
  $\gamma_k \equiv \gamma$, then solution $\hat{\x}$ output by \textbf{Algorithm \ref{alg:online}} using explicit subgradient update
\eqref{eqn:explicit} satisfies
  \begin{align}
    \mathbb{E} [ \rho ( \hat{\x} ) ] & \leq
    \mathcal{O} ( m \log n + \sqrt{n} \log n ) + \myfrac{m ( \bar{a} + \du
    )^2  n }{2} \cdot \gamma \nonumber\\
    \mathbb{E} [ v ( \hat{\x} ) ] & \leq \myfrac{m
    ( \bar{a} + \du )^2}{\dl} + \sqrt{m} ( \bar{a} + \du
    ) + \myfrac{\bar{c}}{\dl} \cdot \gamma^{-1}, \nonumber
  \end{align}
  where $\mathbb{E}[\cdot]$ is  taken over the random permutation.
  }
\end{lem}

\begin{rem}
  \textbf{Lemma \ref{lem:explicit-bound}} implies a trade-off between $\rho$
  and $v$ with respect to the stepsize of the subgradient update. Since we know
  $\bar{a}, \bar{c}, \dl, \du$ in the offline case, it is, therefore, possible
  to find an optimal $\gamma$ to balance the two criteria.
\end{rem}

\begin{thm} \label{thm:explicit-theorem}
  {\tmem{Under the same conditions as \textbf{Lemma {\ref{lem:explicit-bound}}}, if we take
  $\gamma = \sqrt{\myfrac{2 \bar{c}}{\dl ( \bar{a} + \du )^2 m n}}$,
  then \textbf{Algorithm \ref{alg:online}} using explicit subgradient update
\eqref{eqn:explicit} outputs a solution $\hat{\x}$ satisfying
\begin{align*}
	\mathbb{E} [ \rho ( \hat{\x} ) ] &\leq
    \mathcal{O} ( m \log n + \sqrt{n} \log n ) + \Big( \myfrac{( \bar{a} +
    \du )^2 \bar{c}}{2 \dl} \Big)^{1 / 2} \sqrt{m n} \\
    \mathbb{E} [ v ( \hat{\x} ) ] &\leq \myfrac{m
    ( \bar{a} + \du )^2}{\dl} + \sqrt{m} ( \bar{a} + \du
    ) + \Big( \myfrac{( \bar{a} + \du )^2 \bar{c}}{2 \dl}
    \Big)^{1 / 2} \sqrt{m n},
\end{align*}
}
    where $\mathbb{E}[\cdot]$ is  taken over the random permutation.}
\end{thm}

\begin{rem}
  We see the online algorithm gives $\mathcal{O} ( m \log n + \sqrt{n} \log n +
  \sqrt{m n} )$ gap and $\mathcal{O} ( m + \sqrt{m n} )$
  violation even if we take $\gamma$ to be a suboptimal value of $\mathcal{O}
  (1/\sqrt{m n})$. Compared with the ${\mathcal{O}} (\sqrt{n} \log n + m\sqrt{n})$ gap and $\mathcal{O} ( m \sqrt{n} )$ violation of
  {\cite{10.48550/arxiv.2003.02513}}, the result is improved with respect to $m$.
  The main reason why {\cite{10.48550/arxiv.2003.02513}} gets a suboptimal result
  is that they choose $\gamma = 1/\sqrt{n}$, which
  gives an unbalanced trade-off between gap and violation.
\end{rem}

\begin{rem}
The online explicit updates can be
  implemented in $\mathcal{O} ( \text{nnz} ( \A ) )$
  time and is free of any matrix-vector operations. See \textbf{Section \ref{app:impl}} for more details.
\end{rem}

\subsection{Online Implicit Update}
\label{sec:theory-implicit}
In this section, we analyze the performance of the online implicit update
applied to offline LPs. Instead of approximating $f(\y, k)$, we keep all the information using 
\begin{align}
f_{\y^k} ( \y, k ) = f ( \y, k ) = \langle
   \tmd, \y \rangle + [ c_k - \langle \tma_k, \y \rangle
   ]_+ 
	\label{eqn:implicit}
\end{align}
and the primal solution is estimated through proximal point
\begin{eqnarray*}
  \min_{\y, s} & \langle \tmd, \y \rangle + s + \myfrac{1}{2\gamma_k}
  \| \y - \y^k \|^2 & \\
  \text{subject to} & \y \geq \0, s \geq c_k - \langle \tma_k, \y
  \rangle, s \geq 0 & 
\end{eqnarray*}
and we let $x^k = \lambda(s \geq c_k - \langle
\tma_k, \y \rangle)$ be the Lagrangian multiplier of $s \geq c_k - \langle
\tma_k, \y \rangle$ in the optimal solution.

\begin{rem}
  In the literature of online algorithms, the proximal point update is known as
  implicit update, where ``implicit'' comes from the optimality condition
  \[ \0 \in \partial f ( {\y^{k + 1}} , k ) + \gamma_k^{-1} (
     \y^{k + 1} - \y^k ) +\mathcal{N}_{\mathbb{R}_+^n} ( \y^{k + 1}
     ) \]
  and $\y^{k + 1}$ can be expressed implicitly as $$\y^{k + 1} \in \y^k -
  \gamma_k(\partial f ( {\y^{k + 1}} , k )
  +\mathcal{N}_{\mathbb{R}_+^n} ( \y^{k + 1} )),$$ where
  $\mathcal{N}_{\mathbb{R}_+^n} ( \y^{k + 1} )$ is the normal cone
  of the nonnegative orthant at $\y^{k + 1}$. Unlike explicit update which linearizes $f (\y, k)$, implicit update preserves all the information and it is shown to be more robust to stepsize selection
  than subgradient {\cite{10.48550/arxiv.1810.05633, deng2021minibatch}}.
\end{rem}
  
The analysis of implicit update is similar, and we still have a trade-off between gap and violation
\begin{lem} \label{lem:implicit-bound}
  {{Under assumptions {\ref{A1}} to \ref{A3}, if we take
  $\gamma_k \equiv \gamma$, then solution $\hat{\x}$ output by \textbf{Algorithm \ref{alg:online}} using implicit proximal point update
\eqref{eqn:implicit} satisfies
  \begin{align}
    \mathbb{E} [ \rho ( \hat{\x} ) ] & \leq
    \mathcal{O} ( m \log n + \sqrt{n} \log n ) + \myfrac{5m ( \bar{a} + \du
    )^2  n}{2} \cdot \gamma \nonumber\\
    \mathbb{E} [ v ( \hat{\x} ) ] & \leq \myfrac{3m
    ( \bar{a} + \du )^2}{\dl} + \sqrt{m} ( \bar{a} + \du
    ) + \myfrac{\bar{c}}{\dl} \cdot \gamma^{-1}, \nonumber
  \end{align}
  where $\mathbb{E}[\cdot]$ is  taken over the random permutation.}}
\end{lem}

Choosing $\gamma$ properly, we get bounds on gap and violation.

\begin{thm} \label{thm:implicit-theorem}
  {{\tmem{Under the same condition as \textbf{Lemma {\ref{lem:implicit-bound}}}}, if we take
  $\gamma = \sqrt{\myfrac{2 \bar{c}}{5\dl ( \bar{a} + \du )^2 m n}}$,
  then \textbf{Algorithm \ref{alg:online}} using implicit proximal point update
\eqref{eqn:implicit} outputs a solution $\hat{\x}$ satisfying
\begin{align*}
	\mathbb{E} [ \rho ( \hat{\x} ) ] & \leq
    \mathcal{O} ( m \log n + \sqrt{n} \log n ) + \Big( \myfrac{5( \bar{a} +
    \du )^2 \bar{c}}{2 \dl} \Big)^{1 / 2} \sqrt{m n}\\
    \mathbb{E} [ v ( \hat{\x} ) ] & \leq \myfrac{3m
    ( \bar{a} + \du )^2}{\dl} + \sqrt{m} ( \bar{a} + \du
    ) + \Big( \myfrac{5( \bar{a} + \du )^2 \bar{c}}{2 \dl}
    \Big)^{1 / 2} \sqrt{m n},
\end{align*}
    where $\mathbb{E}[\cdot]$ is taken over the random permutation.}}
\end{thm}

\begin{rem}
	Although we do not see an improvement of bounds using implicit update, as will be shown by our numerical experiments, the implicit update often behaves better empirically. One intuitive explanation is that, unlike the explicit update which only outputs binary values, the implicit update is capable of dealing with fractional values. To illustrate this consider the following LP
	\begin{align}
		\max_{0 \leq x_1, x_2 \leq 1}~  x_1 + x_2 ~~ \text{subject to ~} & x_1 + x_2 \leq 0.5 \nonumber
	\end{align}
whose optimal value is $0.5$. If we do not allow constraint violation greater than 0.1, then subgradient update will never take $x=1$ from arbitrary start. However, for the implicit update, we recover the optimal solution if $\gamma < 0.01$.
\end{rem}

\subsection{Improvement by Variable Duplication}
\label{sec:theory-duplicate}
Despite the efficiency of online algorithms, the current optimality gap and constraint violation guarantees are far from enough to solve offline LPs, even approximately. To address this issue, we take advantage of the offline setting and
propose to make a trade-off between time and accuracy, more specifically, by \textbf{1)}. duplicating each column $K$ times \textbf{2)}. running online algorithm on the augmented problem with $nK$ variables \textbf{3)}. taking the average of
primal estimates. This scheme ends up giving an intuitive $\sqrt{K}$ reduction in the final bound. \condsep


\begin{algorithm} 
\caption{Online algorithm with duplication\label{alg:online-dup}}	
\KwIn{$\y^0, K, \{f(\y, j)\}$ from LP data $\A, \tmb, \tmc$}
\begin{itemize}[leftmargin=10pt]
\setlength\itemsep{0em}
	\item Duplicate each of $f(\y, j)$ $K$ times and generate permutation for $nK$ columns (including $\tmc$)
	$\{\underbrace{\tma_{1,1}, \ldots, \tma_{1, n}}_{\text{Duplication 1}}, \underbrace{\tma_{2,1}, \ldots, \tma_{2,n}}_{\text{Duplication 2}}, \ldots, \underbrace{\tma_{K,1}, \ldots, \tma_{K,n}}_{\text{Duplication $K$}}\}$
	\item Run \textbf{Algorithm \ref{alg:online}} and get $\hat{\x}_{\text{Dup}} = (\hat{\x}_1, \ldots, \hat{\x}_K)$
	\item Aggregate solution $\hat{\x} = \frac{1}{K}\sum_{k=1}^K \hat{\x}_k$.
\end{itemize}
\KwOut{$\hat{\x}$}
\end{algorithm}

\begin{thm} \label{thm:dup} 
  Under assumptions \ref{A1} to \ref{A3}, if we apply \textbf{Algorithm \ref{alg:online-dup}} with $K$ duplications and take $\gamma=\mathcal{O}(1/\sqrt{Kmn})$, then both explicit and implicit updates output solutions $\hat{\x}$ satisfying
  \begin{align*}
    \mathbb{E} [ \rho ( \hat{\x} ) ] & \leq
    \mathcal{O} \big( \myfrac{m \log n}{K}+ \sqrt{\myfrac{n}{K}} \log n + \sqrt{\myfrac{m n}{K}}\big) \\
    \mathbb{E} [ v ( \hat{\x} ) ] & \leq \mathcal{O}\big(\myfrac{m}{K} + \sqrt{\myfrac{m n}{K}} \big).
  \end{align*}
\end{thm}

\begin{rem}
There is an extra $\mathcal{O}( \frac{\log K}{\sqrt{K}})$ term in the bound for $\rho$, but note that we generally take $K = \mathcal{O}(n)$ and we drop this term when presenting our results.
\end{rem}

\begin{rem}
We provide an intuitive explanation of the improvement. Taking the explicit update as an example: when each variable is duplicated $K > 1$ times, the final output $\hat{x}^k$ would be allowed to take ${i}/{K}$ for $i \leq K$, while $K = 1$ only allows $\hat{x}^k \in \{0, 1\}$. In other words, larger $K$ offers higher granularity to approximate fractional solutions.
\end{rem}

Till now we have presented the theoretical results of the online algorithms applied to offline LPs. In the following sections, we will focus on how the online algorithms help exact LP solving through sifting. 

\section{Application: Sifting for Linear Programs} \label{sec:sifting}

In the previous sections, we have discussed the use of online algorithms to
approximately solve offline LPs. However, unlike the simplex method, first-order methods rarely
give accurate basis status and generally cannot
be applied for exact LP solving. In this section, we try to alleviate this
issue by identifying the use of our method in sifting, a ``column generation'' framework for linear programming. 
For ease of exposition, we temporarily switch to standard-form LPs $$\max_{\x}~~\langle \tmc, \x \rangle \quad \text{subject to~~} \A \x = \tmb, \x \geq \0.$$
\paragraph{Sifting for LP}
Sifting initially appeared in \cite{forrest1989mathematical} and was formally presented in {\cite{bixby1992very}} to solve LPs with $n \gg m$. Mature LP solvers often use sifting as a candidate solver \cite{pedroso2011optimization, manual1987ibm, ge2022cardinal}. Using the idea
that the set of basic columns $\mathcal{B}= \{ j : x^{\ast}_j > 0 \}$ is small
relative to $[n]$, sifting solves a sequence of working problems which restrict the
problem to a more tractable subset of columns $\mathcal{W} \subseteq [n]$, $| \mathcal{W} |
\ll n$
$$\max_{\x_{\mathcal{W}}}~~\langle \tmc, \x_{\mathcal{W}} \rangle \quad \text{subject to~~} \A \x_{\mathcal{W}}  = \tmb, \x_{\mathcal{W}}  \geq \0.$$
Let $\y_{\mathcal{W}}^{\ast}$
denote the optimal dual solution to the working problem. If
$\y_{\mathcal{W}}^{\ast}$ is dual feasible for the original LP, say,
$\A^{\top} \y_{\mathcal{W}}^{\ast} \geq \tmc$, then $\mathcal{B}
\subseteq \mathcal{W}$ (assuming the optimal solution is unique) and the original problem is solved. Otherwise we price
out the dual infeasible columns $\mathcal{I}= \left\{ j : \left\langle \tma_j,
\y_{\mathcal{W}}^{\ast} \right\rangle < c_j \right\}$ and add them to
$\mathcal{W}$. 

\begin{algorithm} \label{alg:sifting}
\caption{Sifting procedure for LPs}	
\KwIn{Initial working set $\mathcal{W}$}
\While{$\mathcal{I}$ \rm{is nonempty}}{
Solve working problem 

$\max_{\A_{\mathcal{W}} \x_{\mathcal{W}} = \tmb, \x_{\mathcal{W}} \geq \0} 
   ~\langle \tmc_{\mathcal{W}}, \x_{\mathcal{W}} \rangle$
and get $\y_{\mathcal{W}}^{\ast}$

Update $\mathcal{I}= \left\{ j : \left\langle \tma_j,
\y_{\mathcal{W}}^{\ast} \right\rangle < c_j \right\}$
   
$\mathcal{W}=\mathcal{W} \cup \mathcal{I}$
}
\KwOut{Optimal solution to LP}
\end{algorithm}

%

As a special case of column generation, sifting faces challenges that are similar to
column generation {\cite{lubbecke2010column}}
\begin{itemize}
  \item (heading-in) a good initialization of $\mathcal{W}$ is often hard
  \item (dual-oscillation) solution $\y_{\mathcal{W}}^{\ast}$ is unstable
  at the end
\end{itemize}

\paragraph{Dual Stabilization}Among the techniques for sifting, dual stabilization \cite{du1999stabilized, amor2004stabilization, pessoa2018automation}  has been one of the most successful attempts. In a word, most of the dual stabilization techniques work by finding some anchor point $\hat{\y}$ that lies at the center of dual feasible region \cite{2011Chebyshev, gondzio2013new} and then by stablilizing the dual iterations by taking convex combination 
\begin{equation}
	\hat{\y}_{\mathcal{W}}=\alpha {\y}_{\mathcal{W}} + (1 - \alpha) \hat{\y}. \label{eqn:dual-stabilize}
\end{equation}
However, computing an interior point or some center of the dual feasible region is often too costly for huge LPs.

\paragraph{Accelerate Sifting by Online Algorithm}
Generally speaking, to accelerate sifting one needs
\begin{itemize}[itemsep=0pt]
	\item a good estimate of $\{ j : x^{\ast}_j > 0 \}$.
	\item some approximate $\hat{\y} \approx \y^{\ast}$
	\item an inexpensive algorithm that obtains them
\end{itemize}
To this point it's not hard to see that our online algorithm is a perfect candidate to fulfill all three requirements above. First, we can use $\hat{\x}$ as a score function to build initial $\mathcal{W}$. For example, we can choose some threshold value $\tau$ and initialize $\mathcal{W} = \{j: \hat{x}^j \geq \tau\}$. During the sifting
procedure, we can use the approximate dual solution from online algorithm to stabilize the dual iterations. Most importantly, online algorithm is effcient and its running time is negligible compared to the whole sifting procedure.\condsep

In practice, we can efficiently embed sifting into the LP pre-solver by going
through the problem data several times before sifting starts. As our experiments suggest,
online algorithm manages to identify the basis status for a wide range of real-life LPs.\condsep

\begin{rem}
Though the online algorithms work for problems with inequality constraints and upper-bounded variables, 
as the experiments suggest, it often suffices to use the online algorithm as a heuristic in practice even if these constraints are not necessarily satisfied.
\end{rem}

Due to space limit we leave a more detailed discussion of sifting to \textbf{Section \ref{app:sifting}}.

\section{Experiments} \label{sec:experiments}

In this section, we conduct numerical experiments to validate the efficiency of our proposed methods. The experiment is divided into two parts. In the first part we
verify our theoretical results on multi-knapsack benchmark dataset and its
variants. In the second part, we turn to exact LP solving and test
large-scale LP benchmark datasets to see how online algorithms benefit sifting solvers.

\subsection{Approximate Solver}

In this part we test online algorithms in approximate LP solving. We also compare 
the performance of explicit and implicit updates in practice.

\paragraph{Data Description} \label{data-generation} We use synthetic data from multi-knapsack benchmark. More detailedly, we generate benchmark multi-knapsack
problems $\max_{\A \x \leq \tmb, \0 \leq \x \leq \1} \left\langle \tmc, \x
\right\rangle$ as discussed in {\cite{chu1998genetic}}: we
generate each element of $a_{i j}$ uniformly from $\{ 1, \ldots, 1000 \}$.
After generating $\A$, we zero out each element of $\A$ with probability $(1 -
\sigma)$, where $\sigma \in (0, 1]$ controls the sparsity of $\A$. $\tmb$ is
generated by $b_i = \myfrac{\tau}{n} \sum_j a_{i j}$, where $\tau$ is called
tightness coefficient. Each
element of $\tmc$ is generated by $c_i = \myfrac{1}{m} \sum_i a_{i j} +
\delta_i$, where $\delta_i$ is sampled uniformly from $\{ 1, \ldots, 500 \}$.

\paragraph{Performance Metric}
In the first part of our experiment, given a feasible approximate solution
$\hat{\x}$, we use relative optimality to measure its quality
\[ r ( \hat{\x} ) \assign \Big| \myfrac{\langle \tmc,
   \hat{\x} \rangle}{\langle \tmc, \x^{\ast} \rangle}
   \Big| . \]
\paragraph{Testing Configuration and Setup}
\begin{enumerate}[itemsep=0pt,label=\textbf{\arabic*).},ref=\rm{\textbf{A\arabic*}},leftmargin=*,]
  \item \tmtextbf{Dataset}. We  test  $(m, n) \in \{ (5, 10^2), (8, 10^3), (16,
  2\times10^3), (32,   4\times10^3) \}, \tau \in [10^{-2}, 1], \sigma = 1$.
  \item \tmtextbf{Initial Point}. We let online algorithms start from $\0$.
  \item \tmtextbf{Feasibility}. We force the algorithms to respect constraints (take $\hat{x}^k = 1$ only if $\tmb - \sum_{j = 1}^k \tma_j x_j
  \geq \0$).
  \item \tmtextbf{Duplication}. We allow $K \in \{ 1, 2, 4, 8, 16, 32 \}$
  \item \tmtextbf{Stepsize}. We take $\gamma = 1 / \sqrt{K m n}$.
  \item \tmtextbf{Subproblem}. We discuss the way to efficiently compute implicit update in the appendix \textbf{Section \ref{app:subproblem}}.
\end{enumerate}

\begin{figure*}[h]
\centering
\includegraphics[scale=0.22]{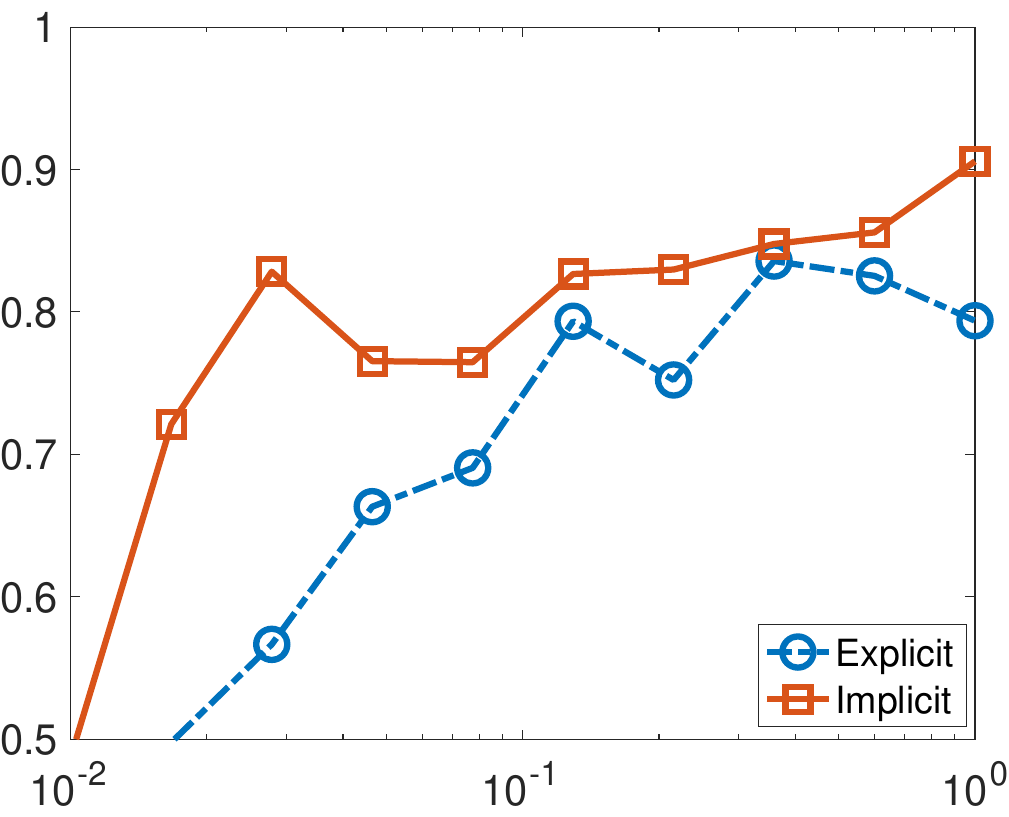}
\includegraphics[scale=0.22]{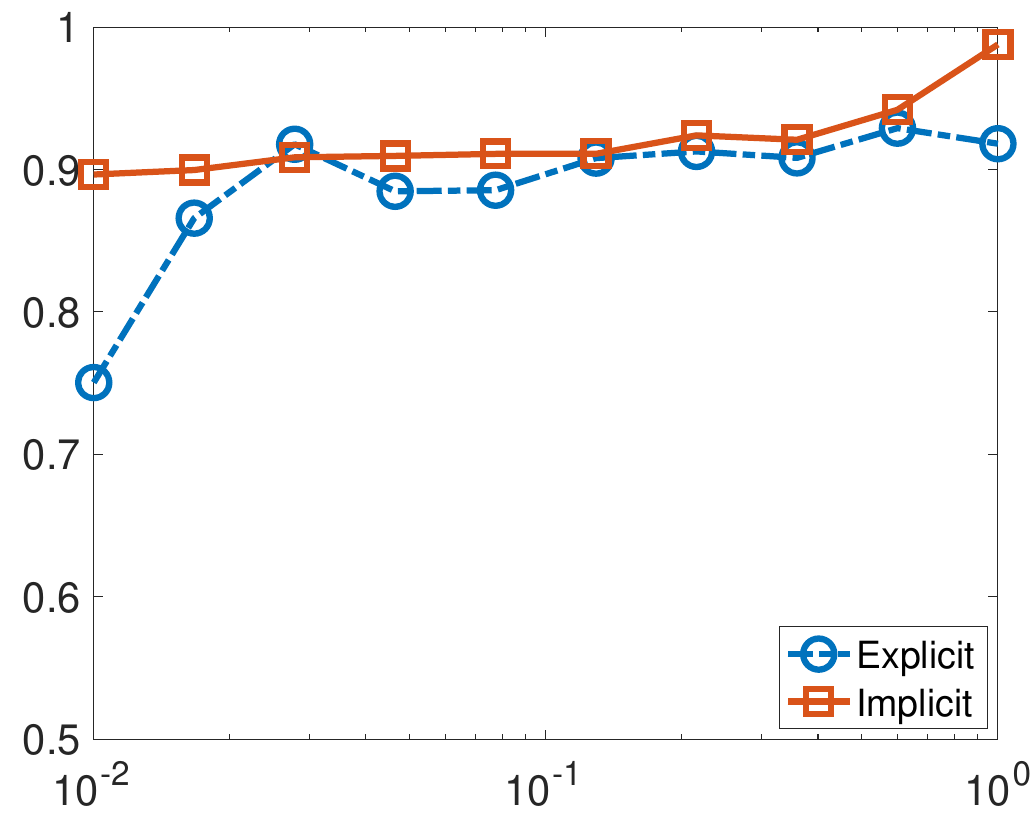}
\includegraphics[scale=0.22]{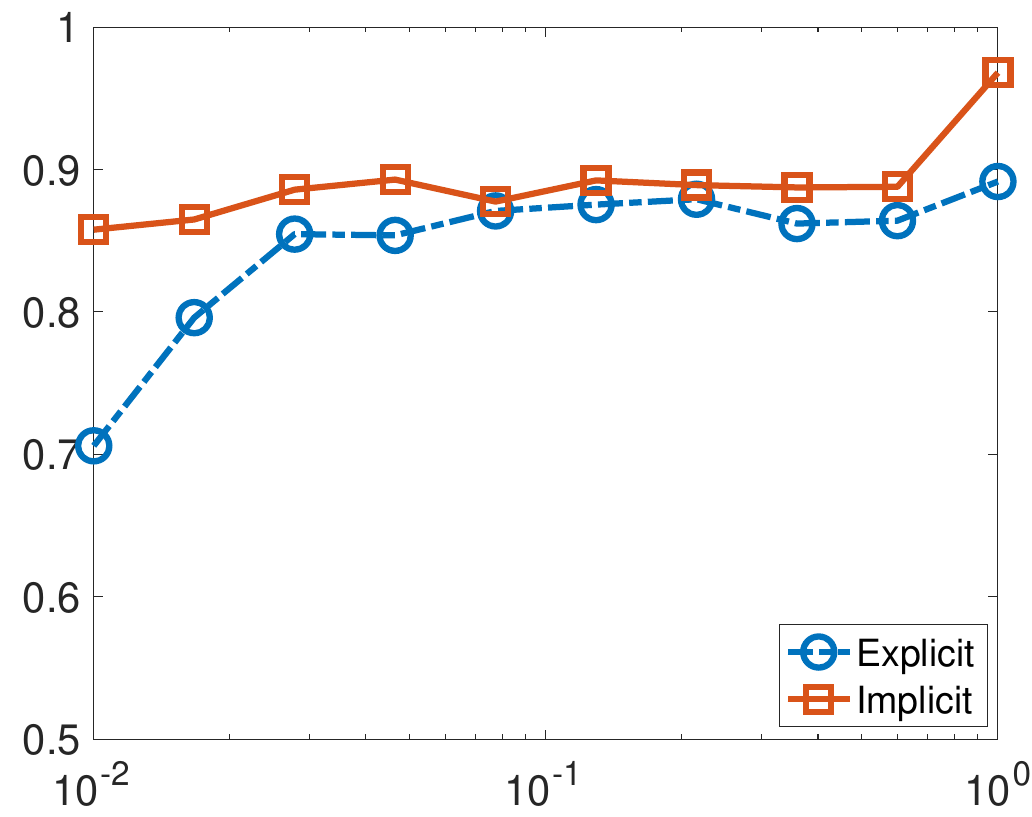}
\includegraphics[scale=0.22]{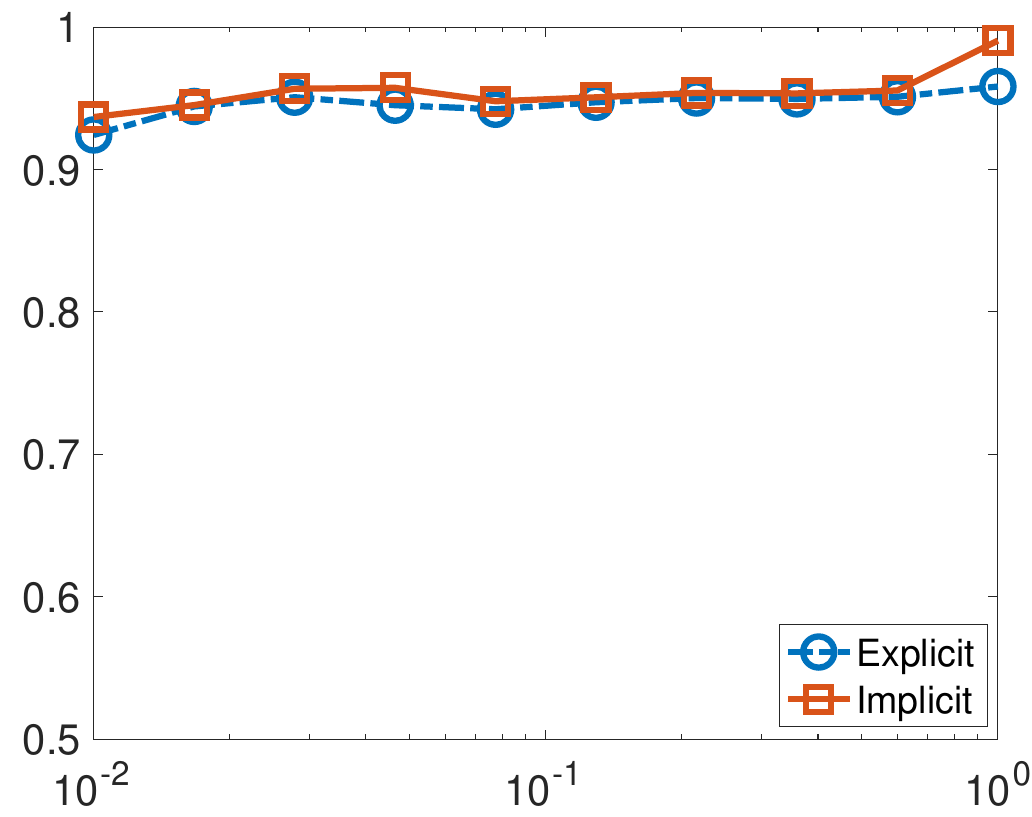}
\includegraphics[scale=0.22]{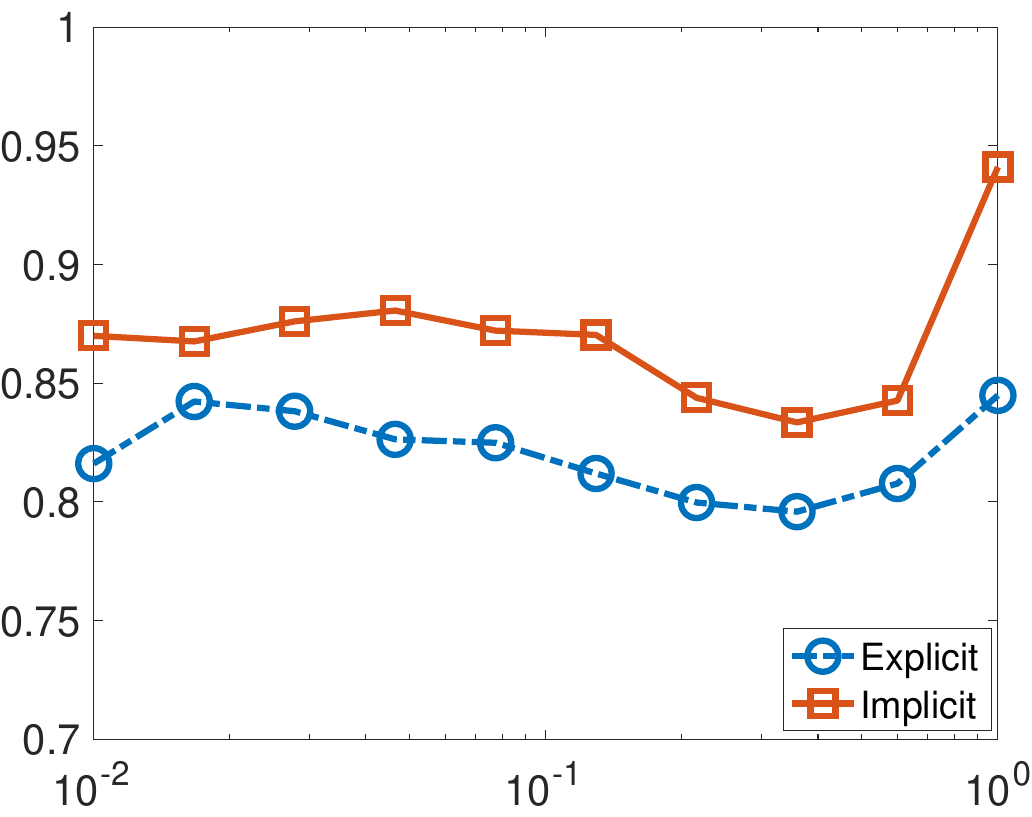}
\includegraphics[scale=0.22]{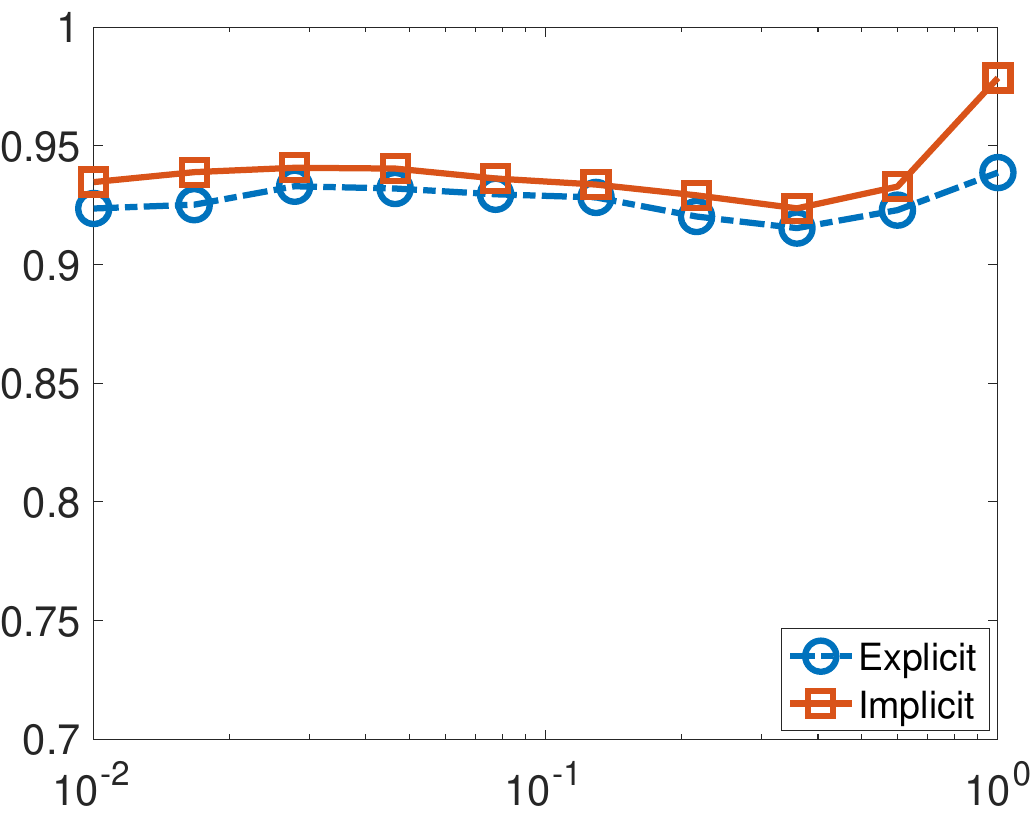}
\includegraphics[scale=0.22]{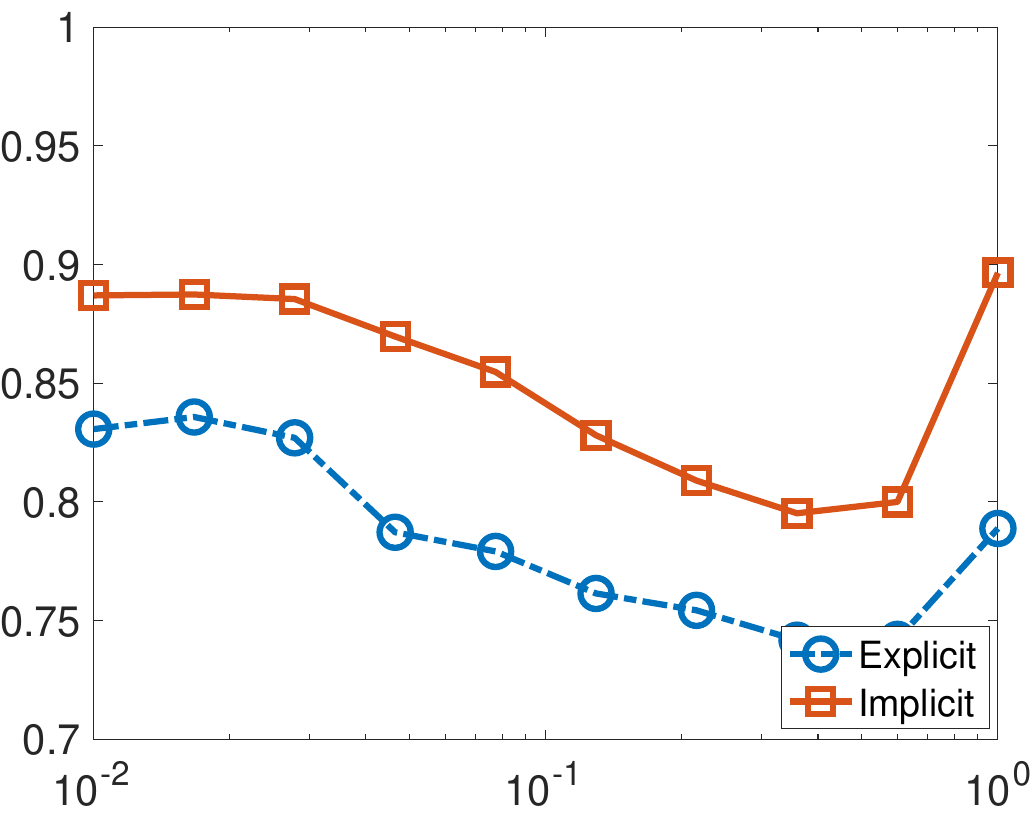}
\includegraphics[scale=0.22]{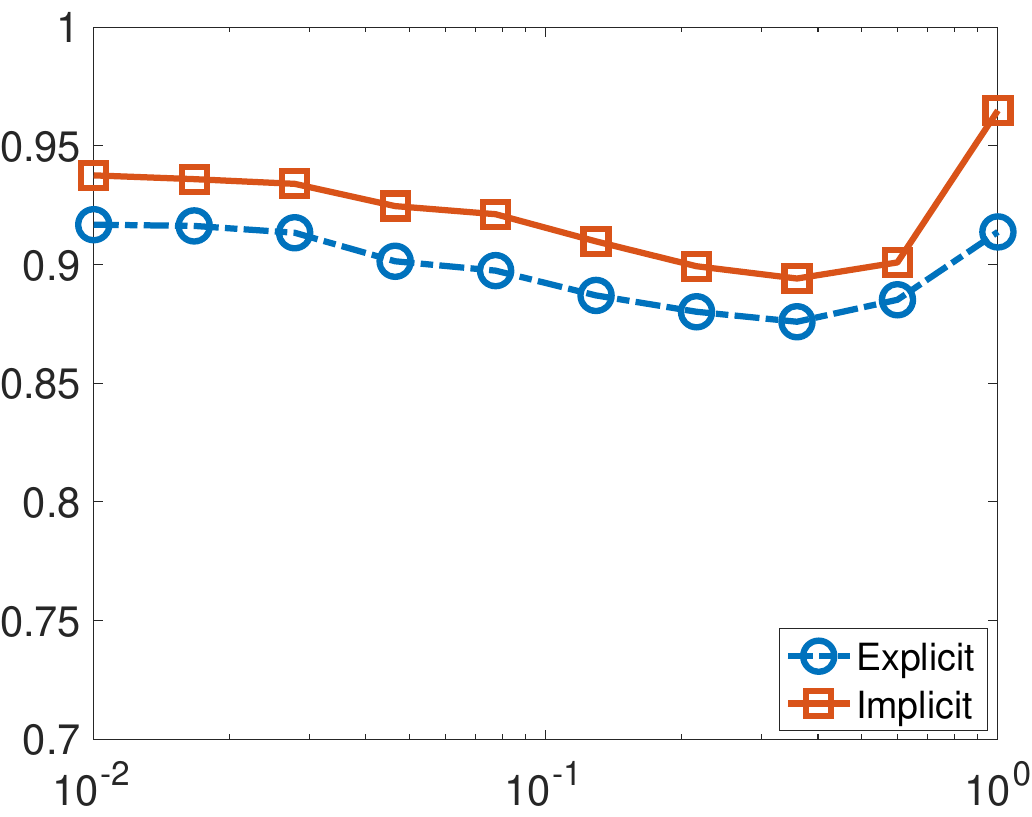}
\caption{\label{fig:1}First row from left to right $(m, n, K) \in \{(5, 10^2, 1), (5, 10^2, 8), (8, 10^3, 1), (8, 10^3, 8)\}$.
Second row from left to right $(m, n, K) \in \{(16, 2\times10^3, 1), (16, 2\times10^3, 8), (32, 4\times10^3, 1), (32, 4\times10^3, 8)\}$.
The x-axis represents $\tau$ parameter ranging from $10^{-2}$ to $1$; The y-axis represents the relative optimality.}
\end{figure*}

First we compare the performance of the two algorithms through their relative optimality under different data settings. For LPs of the same size, we fix $(m, n, K)$ and test $\tau \in [10^{-2}, 1]$ for ten values evenly distributed on the log-scale.  Figure \ref{fig:1} illustrates the performance of both explicit and implicit updates. It can be clearly seen that implicit update, in most cases, outperforms the explicit update. Especially when the constraint is tight ($\tau$ close to 0) and $K = 1$, we see that the performance of implicit updates are dominant. Besides, we observe that with $K$ increasing, explicit update starts to catch up. \condsep

\begin{figure*}[h]
\centering
\includegraphics[scale=0.22]{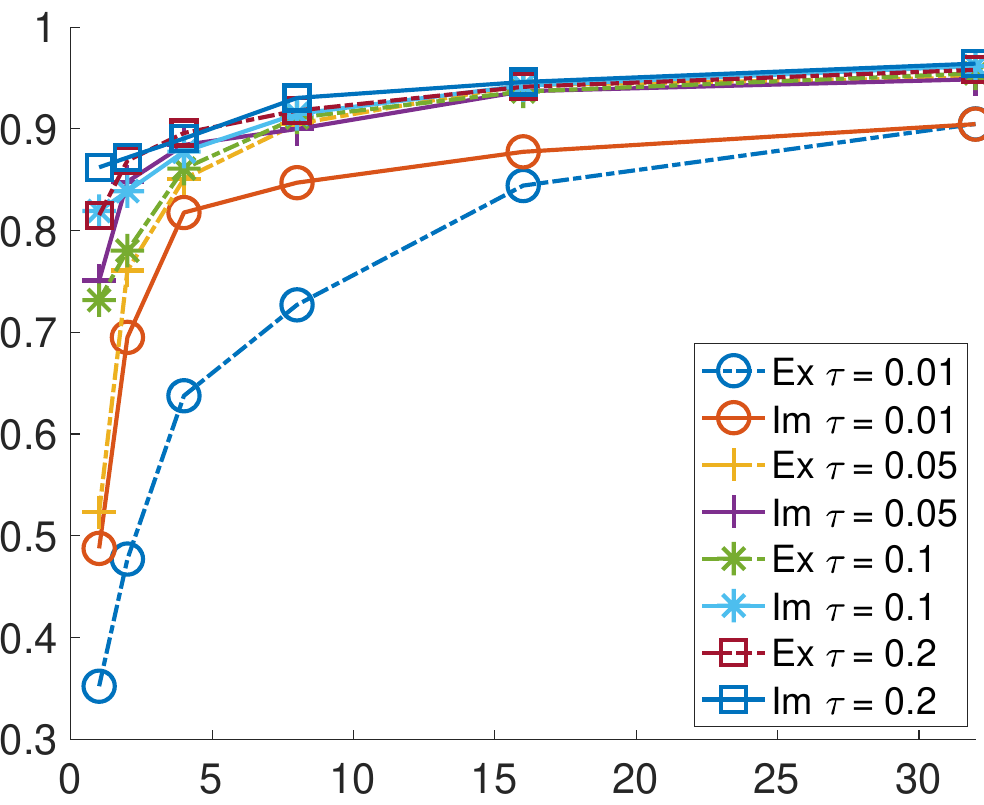}
\includegraphics[scale=0.22]{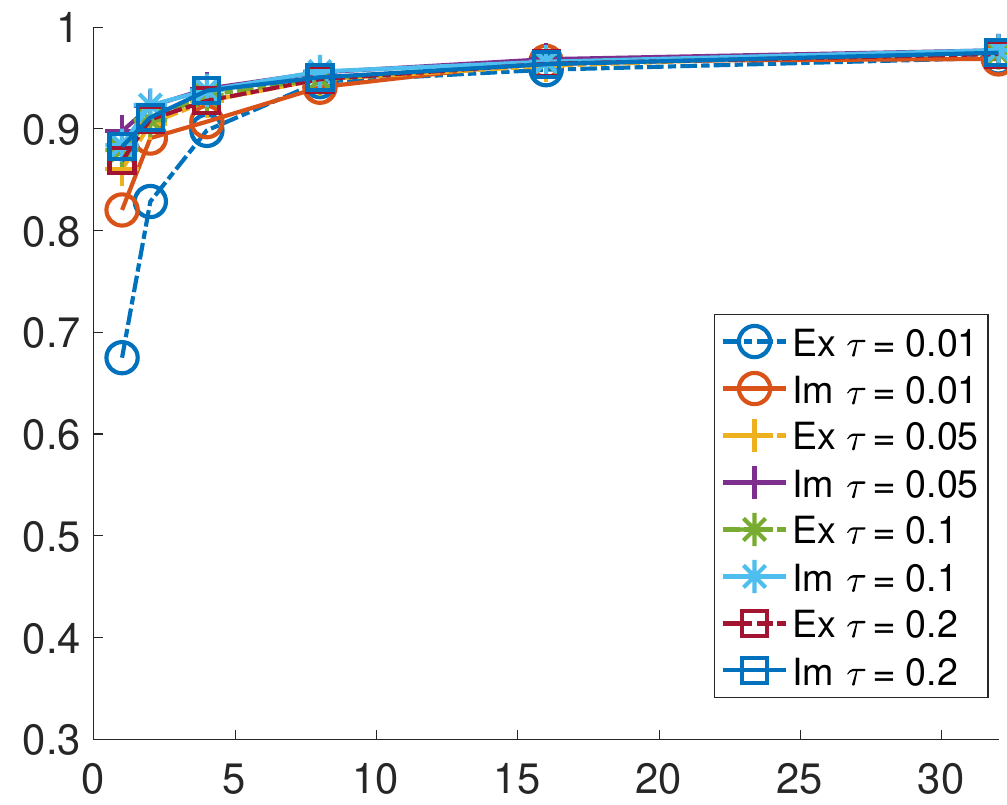}
\includegraphics[scale=0.22]{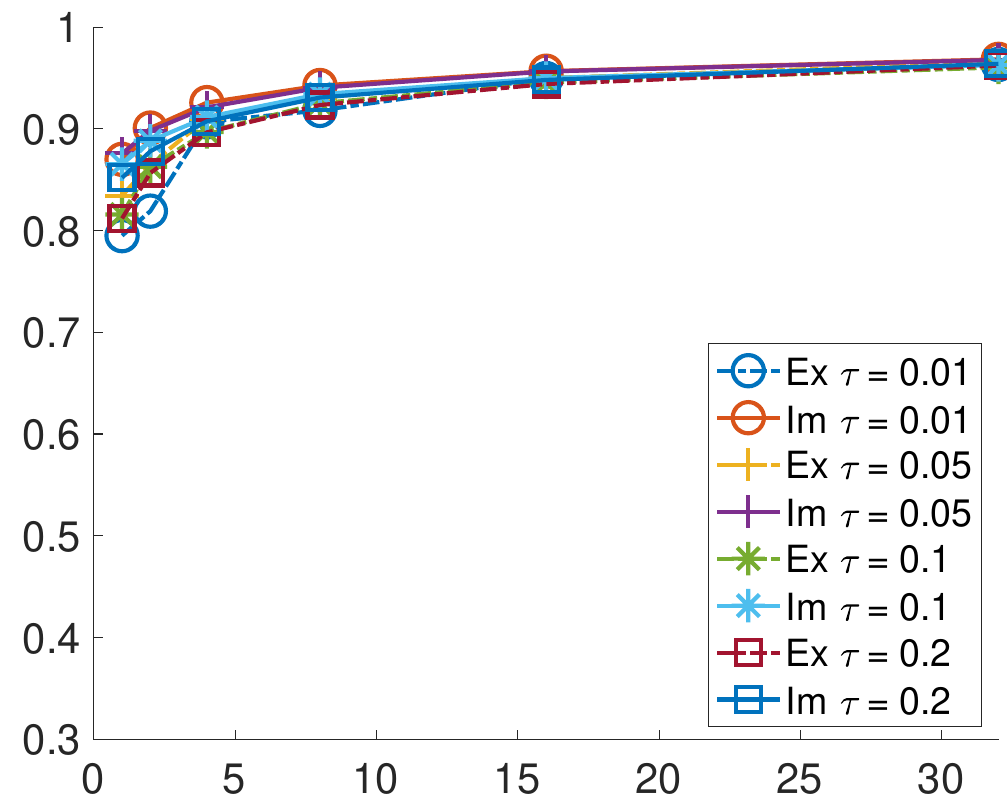}
\includegraphics[scale=0.22]{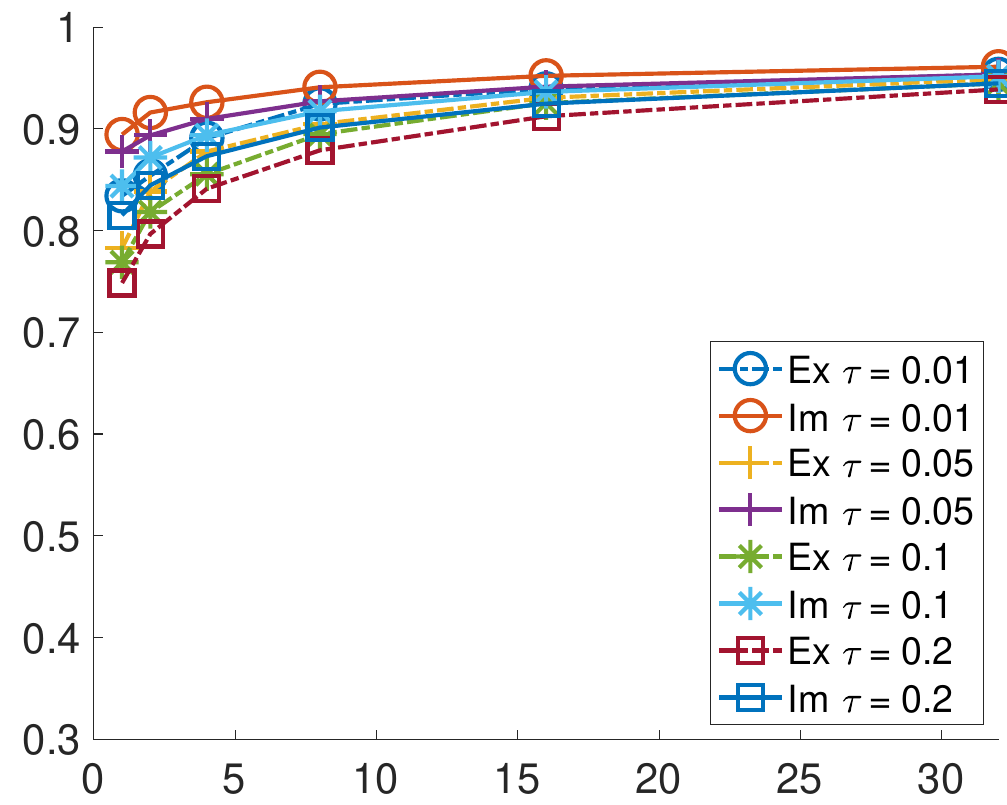}
\caption{\label{fig:2}From left to right $(m, n) \in \{(5, 10^2), (8, 10^3), (16, 2\times10^3), (32, 4\times10^3)\}$. 
The x-axis represents $K$ parameter ranging in $\{1, 2, 4, 8, 16, 32\}$. The y-axis represents the relative optimality.
} 
\end{figure*}

Next we examine the efficiency of the variable duplication scheme by fixing $(m, n, \tau)$ and increasing $K\in\{1, 2, 4, 8, 16, 32\}$. Figure \ref{fig:2} illustrates how variable duplication improves performance of the online algorithms and shows its potential in approximate LP solving. It can be seen that as $K$ increases, relative optimality is gradually improved and we can achieve higher than 90\% relative optimality given moderate $K$. Therefore, it suffices to adopt the online algorithm with variable duplication in the applications where an approximately optimal solution is acceptable. In the appendix \textbf{Section \ref{app:direct}} we further investigate the performance of online algorithms for direct LP solving.
\condsep

Another observation from Figure \ref{fig:2} is that when $\tau$ is close to 0, as we just mentioned, explicit update is dominated by the implicit update due to the restriction of constraint violation. However, as we increase $K$, this gap is soon filled. This suggests in practice we can alternatively combine variable duplication with explicit update to achieve comparable performance to the implicit update. Especially when the constraint matrix is sparse, an $\mathcal{O}(\text{nnz}(\A))$ implementation would be fairly competitive.


\subsection{Sifting and Large LPs}
In this part, we turn to sifting in exact LP solving and see how much online algorithms can help.

\paragraph{Data Description}  We test both synthetic data and real-life large-scale LPs.
For the synthetic data, we generate LP instances in the same way as in \textbf{Section \ref{data-generation}. 
}Besides, we collect 13 instances from {\cite{mittelmann2022benchmark}}.

\begin{table}[h]
\setlength{\tabcolsep}{3.8pt}
  \caption{Datasets collected from large-scale LP benchmark}
   \condsmall
\centering
  \begin{tabular}{cccccc}
  \toprule
    Dataset & Row & Column & Dataset & Row & Column\\
\midrule
    \texttt{rail507} & 507 & 6.4e+04 &
    \texttt{rail516} & 516 & 4.8e+04\\
    \texttt{rail582} & 582 & 5.6e+04 &
    \texttt{rail2586} & 2586 & 9.2e+05\\
    \texttt{rail4284} & 4284 & 1.0e+06 &
    \texttt{scpm1} & 5000 & 5.0e+05\\
    \texttt{scpn2} & 5000 & 1.0e+06 & 
    \texttt{scpl4} & 2000 & 2.0e+05\\
    \texttt{scpj4scip} & 1000 & 1.0e+05 & 
    \texttt{scpk4} & 2000 & 1.0e+05\\
    \texttt{s82} & 87878 & 1.7e+06 & 
    \texttt{s100} & 14733 & 3.6e+05\\
    \texttt{s250r10} & 10962 & 2.7e+05 &
    - & - & - \\
\bottomrule
  \end{tabular}
\end{table}

\paragraph{Performance Metric} We mention in \textbf{Section \ref{sec:sifting}} that online algorithm helps sifting in \textbf{1)}. working problem initialization \textbf{2)}. dual stabilization, and we test these two aspects using different criteria. For working problem initialization, assume that we get an initial working problem estimate $\mathcal{W}$ from the online algorithm and denote by $\mathcal{B}$ the set of basic columns, then we use
$$ \text{acc}(\hat{\x}) := \frac{|\mathcal{B} \cap \mathcal{W}|}{|\mathcal{B}|} \quad \text{and} \quad  \text{rdc}(\hat{\x}):=\frac{|\mathcal{W}|}{n} $$
to respectively evaluate \textbf{1)}. how many basic columns are found out. \textbf{2)}. size of the initialized working problem relative to the original LP.  On the other hand, we also need to evaluate how much, overall, online algorithm can accelerate a sifting solver. To this end we directly use CPU seconds $T(\hat{\x})$ as the metric.\\

\paragraph{Testing Configuration and Setup}
\begin{enumerate}[itemsep=-1pt,label=\textbf{\arabic*).},ref=\rm{\textbf{A\arabic*}},leftmargin=*,]
  \item \tmtextbf{Dataset}. 
We use $(m, n) \in \{(10^2, 10^5), (10^2, 10^6)\}$, $\tau \in \{(0.05, 0.1)\}$ and $\sigma \in \{(0.01, 0.1, 0.15, 0.2)\}$
  \item \tmtextbf{Algorithm Selection}. Since the problems are all sparse, only explicit update is used.
  \item \tmtextbf{Initial Point}. We let online algorithm start from $\1$.
  
  \item \tmtextbf{Duplication}. We take $K = 2$ for all the datasets.
  
  \item \tmtextbf{Stepsize}. We take $\gamma = 1 / \sqrt{K m n}$.
  \item \tmtextbf{Basis Prediction}. Given the output of online algorithm $\hat{\x}$, we use $\mathcal{W} = \{j: \hat{x}^j \geq 1/K\}$ to initialize the working problem.
  \item \tmtextbf{Dual Stabilization}. We implement a basic dual stabilization procedure which takes $\alpha = 0.4$ in \eqref{eqn:dual-stabilize}.
  \item \textbf{Sifting Solver}. We adopt the sifting solver in \texttt{CPLEX 12.10} as our benchmark solver.
\end{enumerate}

\begin{table*}[h]
\centering
\caption{\label{table:1}Performance of initializing working problem and CPU solution time of the large-scale LPs}
\begin{tabular}{ccccc}
  \toprule
  Dataset & $\text{acc}(\hat{\x})$ & $\text{rdc}(\hat{\x})$ & $T_\text{CPLEX}$ & $T(\hat{\x})$\\
  \midrule
  \texttt{rail507} & 271/301 & 11862/62171 & 0.50 & 0.72\\
  \texttt{rail516} & 121/138 & 8572/46978 & 0.38 & 0.73\\
  \texttt{rail582} & 325/347 & 12465/54315 & 0.60 & 1.15\\
  \texttt{rail2586} & 1536/1672 & 145373/909940 & 5.44 & 10.65\\
  \texttt{rail4284} & 1951/2042 & 348135/1090526 & 14.13 & 22.64\\
  \texttt{scpm1} & 2754/2754 & 10352/500000 & 19.08 & 8.74\\
  \texttt{scpn2} & 3411/3411 & 20860/1000000 & 51.38& 12.99 \\
  \texttt{scpl4} & 1149/1149 & 5718/200000 & 0.94 & 0.78\\
  \texttt{scpj4scip} & 552/552 & 3635/99947 & 0.47 & 0.39\\
  \texttt{scpk4} & 930/930 & 4077/100000 & 0.60 & 0.58\\
  \texttt{s82} & 1992/3020 & 52383/1687859 & >3600&>3600\\
  \texttt{s100} & 150/487 & 835/364203 & 151.65&38.30\\
  \texttt{s250r100} & 415/747 & 3080/270323 & 23.07 & 21.55\\
  \bottomrule
\end{tabular}
\end{table*}

\begin{table*}[h]
  \caption{\label{table:2}CPU time of sifting on synthetic datasets}
\centering
  \begin{tabular}{cccccccc}
\toprule
  \multicolumn{4}{c|}{$(m, n) = (10^2, 10^5)$} & \multicolumn{4}{c}{$(m, n) = (10^2, 10^6)$}\\
    \midrule
    $\tau$ & $\sigma$ & $T_{\text{CPLEX}}$ & $T ( \hat{\x} )$ &
    $\tau$ & $\sigma$ & $T_{\text{CPLEX}}$ & $T ( \hat{\x} )$\\
    \midrule
    0.05 & 0.1\% & 8.08 & 3.43 & 0.05 & 0.1\% & 82.79 & 38.20\\
    0.05 & 0.5\% & 7.73 & 4.45 & 0.05 & 0.5\% & 117.02 & 73.39\\
    0.05 & 1\% & 7.86 & 2.81 & 0.05 & 1\% & 61.52 & 35.94\\
    0.05 & 5\% & 7.76 & 3.88 & 0.05 & 5\% & 82.77 & 35.84\\
    0.05 & 10\% & 7.82 & 3.52 & 0.05 & 10\% & 74.89 & 37.67\\
    0.05 & 15\% & 7.96 & 6.30 & 0.05 & 15\% & 65.03 & 40.75\\
    0.05 & 20\% & 7.57 & 3.89 & 0.05 & 20\% & 70.21 & 37.13\\
    0.10 & 0.1\% & 9.19 & 5.57 & 0.10 & 0.1\% & 79.35 & 67.30\\
    0.10 & 0.5\% & 7.57 & 2.71 & 0.10 & 0.5\% & 58.34 & 47.50\\
    0.10 & 1\% & 7.63 & 4.48 & 0.10 & 1\% & 54.69 & 42.73\\
    0.10 & 5\% & 9.08 & 6.51 & 0.10 & 5\% & 56.13 & 52.58\\
    0.10 & 10\% & 8.30 & 5.51 & 0.10 & 10\% & 56.96 & 47.71\\
    0.10 & 15\% & 8.93 & 4.33 & 0.10 & 15\% & 52.40 & 53.79\\
    0.10 & 20\% & 7.82 & 5.18 & 0.10 & 20\% & 57.77 & 42.21\\
    \bottomrule
  \end{tabular}
\end{table*}

Table \ref{table:1} describes the practical performance of the the online algorithm applied to practical large-scale LPs. We can see that on the collected benchmark datasets, online algorithm successfully identifies more than 90\% of the basic columns in most cases while restricting the size of initialization less than 20\% of the original LP size. 
Especially for the \texttt{scp} instances, online algorithm identifies all the basic columns with fewer than 5\% of columns. In this case, the original LP is solved after the first sifting iteration and we only need to solve a much smaller LP. For
the overall LP solving time, we observe a clear speedup on 6 out of the 13 instances, neutral performance on 3 instances and slow down on 4 of the instances. We interpret this slow-down as the effect of our preliminary implementation of the sifting solver compared to \texttt{CPLEX}.\condsep

Finally, we experiment on synthetic datasets that satisfy the assumptions \ref{A1} to \ref{A3}. As \textbf{Table \ref{table:2}} shows, sifting solver, combined with our online algorithm, often outperforms the commercial sifting solvers by more than 50\%. This further illustrates the practical efficiency of our proposed method.

\section{Conclusions}
We adapt two fast online algorithms for offline LPs and obtain algorithms that are free of any matrix multiplication or access to the full LP constraint matrix. We theoretically analyze the optimality gap and constraint violation of the two algorithms and propose a variable-duplication scheme to improve their practical performance. In addition, we identify the potential of online algorithms in exact LP solving when combined with sifting, an LP column generation procedure. Our numerical experiments demonstrate the efficiency of online algorithms, both as an approximate direct solver and as an auxiliary routine in sifting. We believe that it is an interesting direction to introduce online algorithms to the context of offline LP solving.

\section{Acknowledgement}
We thank Yuyang Ye, Xiaocheng Li, and the seminar participants at Shanghai University of Finance and Economics for their fruitful discussions and comments. Especially, we thank Dr. Qi Huangfu for proposing an $\mathcal{O}{(\text{nnz}(\A))}$ implementation of the explicit update method and testing it in the state-of-the-art commercial solver \texttt{COPT} \cite{ge2022cardinal}.

Finally, we sincerely appreciate the efforts the Area Chairs and all the reviewers put into the review process.

\section{Disclosure of Funding}
 This research is partially supported by the National Natural Science Foundation of China (NSFC) [Grant NSFC-72150001, 72225009, 11831002].

\renewcommand \thepart{}
\renewcommand \partname{}

\bibliographystyle{abbrvnat}
\bibliography{olp.bib}

\begin{thebibliography}{50}
\providecommand{\natexlab}[1]{#1}
\providecommand{\url}[1]{\texttt{#1}}
\expandafter\ifx\csname urlstyle\endcsname\relax
  \providecommand{\doi}[1]{doi: #1}\else
  \providecommand{\doi}{doi: \begingroup \urlstyle{rm}\Url}\fi

\bibitem[Agrawal et~al.(2014)Agrawal, Wang, and Ye]{agrawal2014dynamic}
S.~Agrawal, Z.~Wang, and Y.~Ye.
\newblock A dynamic near-optimal algorithm for online linear programming.
\newblock \emph{Operations Research}, 62\penalty0 (4):\penalty0 876--890, 2014.

\bibitem[Amor et~al.(2004)Amor, Desrosiers, and
  Frangioni]{amor2004stabilization}
H.~B. Amor, J.~Desrosiers, and A.~Frangioni.
\newblock \emph{Stabilization in column generation}.
\newblock Groupe d'{\'e}tudes et de recherche en analyse des d{\'e}cisions,
  2004.

\bibitem[Applegate et~al.(2021)Applegate, D{\'\i}az, Hinder, Lu, Lubin,
  O'Donoghue, and Schudy]{applegate2021practical}
D.~Applegate, M.~D{\'\i}az, O.~Hinder, H.~Lu, M.~Lubin, B.~O'Donoghue, and
  W.~Schudy.
\newblock Practical large-scale linear programming using primal-dual hybrid
  gradient.
\newblock \emph{Advances in Neural Information Processing Systems},
  34:\penalty0 20243--20257, 2021.

\bibitem[Asi and Duchi(2019)]{10.48550/arxiv.1810.05633}
H.~Asi and J.~C. Duchi.
\newblock Stochastic (approximate) proximal point methods: Convergence,
  optimality, and adaptivity.
\newblock \emph{SIAM Journal on Optimization}, 29\penalty0 (3):\penalty0
  2257--2290, 2019.

\bibitem[Balseiro et~al.(2020)Balseiro, Lu, and Mirrokni]{balseiro2020dual}
S.~Balseiro, H.~Lu, and V.~Mirrokni.
\newblock Dual mirror descent for online allocation problems.
\newblock In \emph{International Conference on Machine Learning}, pages
  613--628. PMLR, 2020.

\bibitem[Balseiro et~al.(2022)Balseiro, Lu, and
  Mirrokni]{10.48550/arxiv.2011.10124}
S.~R. Balseiro, H.~Lu, and V.~Mirrokni.
\newblock The best of many worlds: Dual mirror descent for online allocation
  problems.
\newblock \emph{Operations Research}, 2022.

\bibitem[Basu et~al.(2020)Basu, Ghoting, Mazumder, and Pan]{basu2020eclipse}
K.~Basu, A.~Ghoting, R.~Mazumder, and Y.~Pan.
\newblock Eclipse: An extreme-scale linear program solver for web-applications.
\newblock In \emph{International Conference on Machine Learning}, pages
  704--714. PMLR, 2020.

\bibitem[Beck(2017)]{beck2017first}
A.~Beck.
\newblock \emph{First-order methods in optimization}.
\newblock SIAM, 2017.

\bibitem[Ben~Amor and Val{\'e}rio~de Carvalho(2005)]{ben2005cutting}
H.~Ben~Amor and J.~Val{\'e}rio~de Carvalho.
\newblock Cutting stock problems.
\newblock In \emph{Column generation}, pages 131--161. Springer, 2005.

\bibitem[Berthold et~al.(2018)Berthold, Farmer, Heinz, and
  Perregaard]{berthold2018parallelization}
T.~Berthold, J.~Farmer, S.~Heinz, and M.~Perregaard.
\newblock Parallelization of the fico xpress-optimizer.
\newblock \emph{Optimization Methods and Software}, 33\penalty0 (3):\penalty0
  518--529, 2018.

\bibitem[Bixby et~al.(1992)Bixby, Gregory, Lustig, Marsten, and
  Shanno]{bixby1992very}
R.~E. Bixby, J.~W. Gregory, I.~J. Lustig, R.~E. Marsten, and D.~F. Shanno.
\newblock Very large-scale linear programming: A case study in combining
  interior point and simplex methods.
\newblock \emph{Operations Research}, 40\penalty0 (5):\penalty0 885--897, 1992.

\bibitem[Bray(2019)]{bray2019logarithmic}
R.~L. Bray.
\newblock Logarithmic regret in multisecretary and online linear programming
  problems with continuous valuations.
\newblock \emph{arXiv e-prints}, pages arXiv--1912, 2019.

\bibitem[Chen et~al.(2022)Chen, Cutkosky, and Orabona]{chen2022implicit}
K.~Chen, A.~Cutkosky, and F.~Orabona.
\newblock Implicit parameter-free online learning with truncated linear models.
\newblock In \emph{International Conference on Algorithmic Learning Theory},
  pages 148--175. PMLR, 2022.

\bibitem[Chu and Beasley(1998)]{chu1998genetic}
P.~C. Chu and J.~E. Beasley.
\newblock A genetic algorithm for the multidimensional knapsack problem.
\newblock \emph{Journal of heuristics}, 4\penalty0 (1):\penalty0 63--86, 1998.

\bibitem[Dantzig(2016)]{dantzig2016linear}
G.~Dantzig.
\newblock Linear programming and extensions.
\newblock In \emph{Linear programming and extensions}. Princeton university
  press, 2016.

\bibitem[Deng and Gao(2021)]{deng2021minibatch}
Q.~Deng and W.~Gao.
\newblock Minibatch and momentum model-based methods for stochastic weakly
  convex optimization.
\newblock \emph{Advances in Neural Information Processing Systems},
  34:\penalty0 23115--23127, 2021.

\bibitem[Deng et~al.(2022)Deng, Feng, Gao, Ge, Jiang, Jiang, Liu, Liu, Xue, Ye,
  et~al.]{deng2022new}
Q.~Deng, Q.~Feng, W.~Gao, D.~Ge, B.~Jiang, Y.~Jiang, J.~Liu, T.~Liu, C.~Xue,
  Y.~Ye, et~al.
\newblock New developments of admm-based interior point methods for linear
  programming and conic programming.
\newblock \emph{arXiv preprint arXiv:2209.01793}, 2022.

\bibitem[Desrosiers and L{\"u}bbecke(2005)]{desrosiers2005primer}
J.~Desrosiers and M.~E. L{\"u}bbecke.
\newblock A primer in column generation.
\newblock In \emph{Column generation}, pages 1--32. Springer, 2005.

\bibitem[Du~Merle et~al.(1999)Du~Merle, Villeneuve, Desrosiers, and
  Hansen]{du1999stabilized}
O.~Du~Merle, D.~Villeneuve, J.~Desrosiers, and P.~Hansen.
\newblock Stabilized column generation.
\newblock \emph{Discrete Mathematics}, 194\penalty0 (1-3):\penalty0 229--237,
  1999.

\bibitem[Forrest(1989)]{forrest1989mathematical}
J.~Forrest.
\newblock Mathematical programming with a library of optimization subroutines.
\newblock In \emph{ORSA/TIMS Joint National Meeting, New York}, 1989.

\bibitem[Fountoulakis and Gondzio(2016)]{fountoulakis2016performance}
K.~Fountoulakis and J.~Gondzio.
\newblock Performance of first-and second-order methods for l1-regularized
  least squares problems.
\newblock \emph{Computational Optimization and Applications}, 65\penalty0
  (3):\penalty0 605--635, 2016.

\bibitem[Garstka et~al.(2021)Garstka, Cannon, and Goulart]{garstka2021cosmo}
M.~Garstka, M.~Cannon, and P.~Goulart.
\newblock Cosmo: A conic operator splitting method for convex conic problems.
\newblock \emph{Journal of Optimization Theory and Applications}, 190\penalty0
  (3):\penalty0 779--810, 2021.

\bibitem[Ge et~al.(2022)Ge, Huangfu, Wang, Wu, and Ye]{ge2022cardinal}
D.~Ge, Q.~Huangfu, Z.~Wang, J.~Wu, and Y.~Ye.
\newblock Cardinal optimizer (copt) user guide.
\newblock \emph{arXiv preprint arXiv:2208.14314}, 2022.

\bibitem[Gondzio et~al.(2013)Gondzio, Gonz{\'a}lez-Brevis, and
  Munari]{gondzio2013new}
J.~Gondzio, P.~Gonz{\'a}lez-Brevis, and P.~Munari.
\newblock New developments in the primal--dual column generation technique.
\newblock \emph{European Journal of Operational Research}, 224\penalty0
  (1):\penalty0 41--51, 2013.

\bibitem[Hazan et~al.(2016)]{hazan2016introduction}
E.~Hazan et~al.
\newblock Introduction to online convex optimization.
\newblock \emph{Foundations and Trends{\textregistered} in Optimization},
  2\penalty0 (3-4):\penalty0 157--325, 2016.

\bibitem[Jiang et~al.(2022)Jiang, Ma, and Zhang]{10.48550/arxiv.2210.07996}
J.~Jiang, W.~Ma, and J.~Zhang.
\newblock {Degeneracy is OK: Logarithmic Regret for Network Revenue Management
  with Indiscrete Distributions}.
\newblock \emph{arXiv}, 2022.
\newblock \doi{10.48550/arxiv.2210.07996}.

\bibitem[Kasirzadeh et~al.(2017)Kasirzadeh, Saddoune, and
  Soumis]{kasirzadeh2017airline}
A.~Kasirzadeh, M.~Saddoune, and F.~Soumis.
\newblock Airline crew scheduling: models, algorithms, and data sets.
\newblock \emph{EURO Journal on Transportation and Logistics}, 6\penalty0
  (2):\penalty0 111--137, 2017.

\bibitem[Kulis and Bartlett(2010)]{kulis2010implicit}
B.~Kulis and P.~L. Bartlett.
\newblock Implicit online learning.
\newblock In \emph{Proceedings of the 27th International Conference on Machine
  Learning (ICML-10)}, pages 575--582, 2010.

\bibitem[Lee and Park(2011{\natexlab{a}})]{2011Chebyshev}
C.~Lee and S.~Park.
\newblock Chebyshev center based column generation.
\newblock \emph{Discrete Applied Mathematics}, 159\penalty0 (18):\penalty0
  2251--2265, 2011{\natexlab{a}}.

\bibitem[Lee and Park(2011{\natexlab{b}})]{lee2011chebyshev}
C.~Lee and S.~Park.
\newblock Chebyshev center based column generation.
\newblock \emph{Discrete Applied Mathematics}, 159\penalty0 (18):\penalty0
  2251--2265, 2011{\natexlab{b}}.

\bibitem[Li and Ye(2022)]{10.48550/arxiv.1909.05499}
X.~Li and Y.~Ye.
\newblock Online linear programming: Dual convergence, new algorithms, and
  regret bounds.
\newblock \emph{Operations Research}, 70\penalty0 (5):\penalty0 2948--2966,
  2022.

\bibitem[Li et~al.(2020)Li, Sun, and Ye]{10.48550/arxiv.2003.02513}
X.~Li, C.~Sun, and Y.~Ye.
\newblock Simple and fast algorithm for binary integer and online linear
  programming.
\newblock \emph{Advances in Neural Information Processing Systems},
  33:\penalty0 9412--9421, 2020.

\bibitem[Lin et~al.(2021)Lin, Ma, Ye, and Zhang]{lin2021admm}
T.~Lin, S.~Ma, Y.~Ye, and S.~Zhang.
\newblock An admm-based interior-point method for large-scale linear
  programming.
\newblock \emph{Optimization Methods and Software}, 36\penalty0 (2-3):\penalty0
  389--424, 2021.

\bibitem[L{\"u}bbecke(2010)]{lubbecke2010column}
M.~E. L{\"u}bbecke.
\newblock Column generation.
\newblock \emph{Wiley encyclopedia of operations research and management
  science}, 2010.

\bibitem[L{\"u}bbecke and Desrosiers(2005)]{lubbecke2005selected}
M.~E. L{\"u}bbecke and J.~Desrosiers.
\newblock Selected topics in column generation.
\newblock \emph{Operations research}, 53\penalty0 (6):\penalty0 1007--1023,
  2005.

\bibitem[Luo and Sun(1998)]{luo1998analytic}
Z.-Q. Luo and J.~Sun.
\newblock An analytic center based column generation algorithm for convex
  quadratic feasibility problems.
\newblock \emph{SIAM Journal on Optimization}, 9\penalty0 (1):\penalty0
  217--235, 1998.

\bibitem[Ma et~al.(2022)Ma, Cao, Tsang, and Xia]{10.48550/arxiv.2209.00399}
W.~Ma, Y.~Cao, D.~H.~K. Tsang, and D.~Xia.
\newblock {Optimal Regularized Online Convex Allocation by Adaptive
  Re-Solving}.
\newblock \emph{arXiv}, 2022.
\newblock \doi{10.48550/arxiv.2209.00399}.

\bibitem[Manual(1987)]{manual1987ibm}
C.~U. Manual.
\newblock Ibm ilog cplex optimization studio.
\newblock \emph{Version}, 12\penalty0 (1987-2018):\penalty0 1, 1987.

\bibitem[Megiddo and Chandrasekaran(1989)]{megiddo1989varepsilon}
N.~Megiddo and R.~Chandrasekaran.
\newblock On the $\varepsilon$-perturbation method for avoiding degeneracy.
\newblock \emph{Operations Research Letters}, 8\penalty0 (6):\penalty0
  305--308, 1989.

\bibitem[Mittelmann(2022)]{mittelmann2022benchmark}
H.~Mittelmann.
\newblock Benchmark of barrier lp solvers, 2022.

\bibitem[Moharrami et~al.(2015)Moharrami, Mahini, and
  Cocchetti]{moharrami2015elastoplastic}
H.~Moharrami, M.~Mahini, and G.~Cocchetti.
\newblock Elastoplastic analysis of plane stress/strain structures via
  restricted basis linear programming.
\newblock \emph{Computers \& Structures}, 146:\penalty0 1--11, 2015.

\bibitem[Nesterov(2012)]{nesterov2012efficiency}
Y.~Nesterov.
\newblock Efficiency of coordinate descent methods on huge-scale optimization
  problems.
\newblock \emph{SIAM Journal on Optimization}, 22\penalty0 (2):\penalty0
  341--362, 2012.

\bibitem[Nesterov(2014)]{nesterov2014subgradient}
Y.~Nesterov.
\newblock Subgradient methods for huge-scale optimization problems.
\newblock \emph{Mathematical Programming}, 146\penalty0 (1):\penalty0 275--297,
  2014.

\bibitem[O'donoghue et~al.(2016)O'donoghue, Chu, Parikh, and Boyd]{o2016conic}
B.~O'donoghue, E.~Chu, N.~Parikh, and S.~Boyd.
\newblock Conic optimization via operator splitting and homogeneous self-dual
  embedding.
\newblock \emph{Journal of Optimization Theory and Applications}, 169\penalty0
  (3):\penalty0 1042--1068, 2016.

\bibitem[Pedroso(2011)]{pedroso2011optimization}
J.~P. Pedroso.
\newblock Optimization with gurobi and python.
\newblock \emph{INESC Porto and Universidade do Porto,, Porto, Portugal}, 1,
  2011.

\bibitem[Perez et~al.(2022)Perez, Ament, Gomes, and
  Barlaud]{perez2022efficient}
G.~Perez, S.~Ament, C.~Gomes, and M.~Barlaud.
\newblock Efficient projection algorithms onto the weighted l1 ball.
\newblock \emph{Artificial Intelligence}, 306:\penalty0 103683, 2022.

\bibitem[Pessoa et~al.(2018)Pessoa, Sadykov, Uchoa, and
  Vanderbeck]{pessoa2018automation}
A.~Pessoa, R.~Sadykov, E.~Uchoa, and F.~Vanderbeck.
\newblock Automation and combination of linear-programming based stabilization
  techniques in column generation.
\newblock \emph{INFORMS Journal on Computing}, 30\penalty0 (2):\penalty0
  339--360, 2018.

\bibitem[Talluri et~al.(2004)Talluri, Van~Ryzin, and
  Van~Ryzin]{talluri2004theory}
K.~T. Talluri, G.~Van~Ryzin, and G.~Van~Ryzin.
\newblock \emph{The theory and practice of revenue management}, volume~1.
\newblock Springer, 2004.

\bibitem[Yanover et~al.(2006)Yanover, Meltzer, Weiss, Bennett, and
  Parrado-Hern{\'a}ndez]{yanover2006linear}
C.~Yanover, T.~Meltzer, Y.~Weiss, K.~P. Bennett, and E.~Parrado-Hern{\'a}ndez.
\newblock Linear programming relaxations and belief propagation--an empirical
  study.
\newblock \emph{Journal of Machine Learning Research}, 7\penalty0 (9), 2006.

\bibitem[Ye(2011)]{ye2011interior}
Y.~Ye.
\newblock \emph{Interior point algorithms: theory and analysis}.
\newblock John Wiley \& Sons, 2011.

\end{thebibliography}

\doparttoc
\faketableofcontents
\part{}

\appendix
\onecolumn
\addcontentsline{toc}{section}{Appendix}
\part{Appendix} 
\parttoc

\vspace{10pt}

\textbf{Structure of the Appendix}

The appendix is organized as follows.
\textbf{Section \ref{app:theory}} introduces auxiliary results and proves our main results from \textbf{Section \ref{sec:theory}}.
As a complement to our discussions in \textbf{Section \ref{sec:experiments}}, \textbf{Section \ref{app:applications}} gives a more comprehensive treatment of online linear programming applied to \textbf{1)}. Direct LP solving. \textbf{2)}. LP sifting and column generation. \textbf{3)}. Integer programming. Particularly we discuss the practical implementation of both explicit and implicit update in \textbf{Section \ref{app:impl}} and \textbf{Section \ref{app:subproblem}}.

\newpage

\section{Theoretical Aspects of Online Linear Programming} \label{app:theory}

\subsection{Auxiliary Results} \label{app:aux}

In this section, we present some auxiliary results from
{\cite{10.48550/arxiv.2003.02513}} that will help in the proof. Recall that given an index set $\mathcal{S} \subseteq [n],$ we use $\A_{\mathcal{S}}$ to denote the sub-matrix indexed from columns of $\A$ and use $\tmc_{\mathcal{S}}$ to denote a sub-vector indexed from $\tmc$. Then we introduce two auxiliary LPs as follows.

\textbf{Auxiliary LPs from Proposition 1 of \cite{10.48550/arxiv.2003.02513}}
\begin{eqnarray*}
  F_s^{\ast} \assign \max_{\x_{[1 : s]}} & \langle \tmc_{[1 : s]}, \x_{[1 : s]}
  \rangle & \\
  \text{subject to} & \A_{[1 : s]} \x_{[1 : s]} \leq \frac{s \tmb}{n} & \\
  & \0 \leq \x_{[1 : s]} \leq \1 & 
\end{eqnarray*}
\begin{eqnarray*}
  \hat{F}_k^{\ast} \assign \max_{\x_{[k : n]}} & \langle \tmc_{[k : n]}, \x_{[k :
  n]} \rangle & \\
  \text{subject to} & \A_{[k : n]} \x_{[k : n]} \leq ( 1 - \tfrac{k -
  1}{n} ) \tmb & \\
  & \0 \leq \x_{[k : n]} \leq \1, & 
\end{eqnarray*}
The optimal values of of the two LPs are denoted by $F_s^{\ast}$ and $\hat{F}_k^{\ast}$
respectively. The following lemma provides a lens to deal with random
permutation.

\begin{lem}[Optimality gap {\cite{10.48550/arxiv.2003.02513}}] \label{lem:li}
  We have the following bound on the optimality gap of the online algorithm
  \begin{align*}
    F^{\ast}_n - \sum_{k = 1}^n \mathbb{E} [c_k x^k] & \leq m \bar{c} +
    \frac{\bar{c}  \sqrt{n} \log n}{\dl} + m \bar{c} \log n + \frac{\alpha
    \bar{c}}{n} + \sum_{k = 1}^n \mathbb{E} \Big[ \tfrac{\hat{F}_{n - k +
    1}^{\ast}}{n - k + 1} - c_k x^k \Big]\\
    & = \mathcal{O} ( m \log n + \sqrt{n} \log n ) + \sum_{k = 1}^n
    \mathbb{E} \Big[ \tfrac{\hat{F}_{n - k + 1}^{\ast}}{n - k + 1} - c_k x^k
    \Big],
  \end{align*}
  where $\alpha = \max \{ e, e^{16 \bar{a}^2}, 16 \bar{a}^2 \}$ and the expectation is taken over the random permutation.
\end{lem}

\begin{rem}
Given \textbf{Lemma \ref{lem:li}}, it remains to analyze the quantity $\sum_{k=1}^n\mathbb{E} \Big[
\frac{\hat{F}_{n - k + 1}^{\ast}}{n - k + 1} - c_k x^k \Big]$ to bound the optimality gap, and in our proof we analyze it under explicit and implicit updates respectively.
\end{rem}

We also need a well-known three-point lemma in the analysis of implicit
updates.

\begin{lem}[Three-point lemma {\cite{beck2017first}}]
  Let $f$ be a convex function and
  \[ \y^+ = \argmin_{\y}  \Big\{ f ( \y) + \frac{1}{2 \gamma}
     \| \y - \mathbf{z} \|^2 \Big\}, \]
  where $\gamma > 0$. Then we have
  \[ f ( \y^+ ) + \frac{1}{2 \gamma} \| \y^+ - \mathbf{z}
     \|^2 \leq f ( \y ) + \frac{1}{2 \gamma} \| \y -
     \mathbf{z} \|^2 - \frac{1}{2 \gamma} \| \y - \y^+ \|^2 .
  \]
\end{lem}

\paragraph{Clue of proof} With the tools in hand, now we get down to the analysis of different algorithms. Our proof basically follows three steps as follows:
\begin{itemize}
  \item Tracking the dual iteration $\{\y^k\}$
  \item Tracking constraint violation and optimality gap
  \item Taking trade-off between optimality gap and constraint violation by properly choosing $\gamma$.
\end{itemize}

\subsection{Proof of Results in Section \ref{sec:theory-explicit}} \label{app:proof-3-1}

In this section, we present the proof for subgradient-based explicit online algorithm.

The following lemma tracks the dual iterations of \textbf{Algorithm \ref{alg:online}} with explicit update.

\begin{lem}[Tracking the dual iteration]
  {\tmem{\label{lem:explicit-dual-iter}
  
Under \ref{A1} and \ref{A2}, if we let $\{ \y^k \}$ be the sequence of dual iterates generated by
\textbf{Algorithm \ref{alg:online}} with explicit update and $\gamma_k \equiv \gamma$, then
  \[ \| \y^k \| \leq \frac{m ( \bar{a} + \du )^2
     \gamma}{\dl} + \sqrt{m} ( \bar{a} + \du
     ) \gamma + \frac{\bar{c}}{\dl}. \]}}
\end{lem}

Boundedness of dual iterations turns out to be important for limiting the constraint violation. Now we state the detailed proof of the theoretical results.

\subsubsection{Proof of Lemma \ref{lem:explicit-dual-iter}}

The proof is adapted from {\cite{10.48550/arxiv.2003.02513}} and is improved
via a sharper analysis. First recall that we update the dual iterations by
\begin{eqnarray*}
  \y^{k + 1}= [ \y^k + \gamma ( \tma_k \mathbb{I} \{ c_k >
  \langle \tma_k, \y^k \rangle \} - \tmd )  ]_+ = [ \y^k + \gamma ( \tma_k x^k - \tmd )  ]_+
\end{eqnarray*}

and we successively deduce that
\begin{align}
  &  \| \y^{k + 1} \|^2 - \| \y^k \|^2 \nonumber \\ 
  & \leq \| \y^k + \gamma ( \tma_k x^k - \tmd ) \|^2
  - \| \y^k \|^2 \label{eqn:lem-5-1} \\ 
  & \leq - 2 \gamma \langle \tmd - \tma_k x^k, \y^k \rangle + m
  ( \bar{a} + \du )^2 \gamma^2 \label{eqn:lem-5-2} \\
  & \leq m ( \bar{a} + \du )^2 \gamma^2 + 2 \gamma \bar{c} - 2
  \gamma \dl \| \y^k \|, \label{eqn:lem-5-3}
\end{align}
where \eqref{eqn:lem-5-1} is due to $\|[\x]_+\| \leq \|\x\|$;
\eqref{eqn:lem-5-2} uses \ref{A1}, \ref{A2} to get $\gamma^2 \| \tma_k x^k - \tmd
\|^2 \leq m ( \bar{a} + \du )^2 \gamma^2$ since $\|
\tma_k \|_{\infty} \leq \bar{a}, \0 < \tmd \leq \du \cdot \1$ and $|x^k| \leq 1$;
\eqref{eqn:lem-5-3} uses the relation $$ - 2 \gamma \langle \tmd, \y^k \rangle
\leq - 2 \gamma \dl \| \y^k \|_1 \leq - 2 \gamma \dl \| \y^k
\|_2$$ and $\langle \tma_k x^k, \y^k \rangle \leq c_k \leq
\bar{c}$. 

On the other hand, we have, by triangle inequality that
\[ \| \y^{k + 1} \| \leq \| \y^k + \gamma ( \tma_k x^k -
   \tmd ) \| \leq \| \y^k \| + \gamma \| \tma_k
   x^k - \tmd \| \leq \| \y^k \| + \gamma \sqrt{m} (
   \bar{a} + \du ) \]
and if $\| \y^k \| \geq \frac{m ( \bar{a}
+ \du )^2 \gamma + 2 \bar{c}}{2 \dl}$,  by \eqref{eqn:lem-5-3} we know that
\[ \| \y^{k + 1} \|^2 - \| \y^k \|^2 \leq m (
   \bar{a} + \du )^2 \gamma^2 + 2 \gamma \bar{c} - 2 \gamma \dl \|
   \y^k \| \leq 0. \]
Therefore, if $\|\y^1\| \leq \frac{m (
\bar{a} + \du )^2 \gamma + 2 \bar{c}}{2 \dl}$, then $\| \y^k \|$ never exceeds $\frac{m (
\bar{a} + \du )^2 \gamma + 2 \bar{c}}{2 \dl} + \gamma \sqrt{m} (
\bar{a} + \du )$ and this completes the proof.

\subsubsection{Proof of Lemma \ref{lem:explicit-bound}}

By the updating formula of subgradient,
\[ \y^{k + 1} = [ \y^k + \gamma ( \tma x^k - \tmd )  ]_+ \geq \y^k + \gamma
   ( \tma_k x^k - \tmd ) \]
and telescoping over $k = 1, \ldots, n$ gives
\[\A \hat{\x} =  \sum_{k = 1}^n \tma_k x^k \leq \tmb + \gamma^{- 1} \sum_{k = 1}^n (
   \y^{k + 1} - \y^k ) \leq \tmb + \gamma^{- 1} \y^{n + 1} . \]
Re-arranging the term, taking positive part and taking expectation with respect to the random permutation, we have
\begin{align}
  \mathbb{E} [ \| [ \A \hat{\x} - \tmb ]_+ \| ] & \leq  \gamma^{- 1} \mathbb{E} [ \| \y^{k + 1} \| ] \nonumber\\
  & \leq  \gamma^{- 1}  \left[ \frac{m ( \bar{a} + \du )^2
  \gamma}{\dl} + \gamma \sqrt{m} ( \bar{a} + \du
  ) + \frac{\bar{c}}{\dl} \right] \label{eqn:lem-1-1} \\
  & = \frac{m ( \bar{a} + \du )^2}{\dl} + \sqrt{m} (
  \bar{a} + \du ) + \frac{\bar{c}}{\gamma \dl}, \nonumber 
\end{align}
where \eqref{eqn:lem-1-1} invokes \textbf{Lemma \ref{lem:explicit-dual-iter}} to bound the constraint violation.

Next we bound $\sum_{k = 1}^n \mathbb{E} \Big[ \frac{\hat{F}_{n - k +
1}^{\ast}}{n - k + 1} - c_k x^k \Big]$ and deduce that
\begin{align}
  & \quad \| \y^{k
  + 1} \|^2 - \| \y^k \|^2 \nonumber\\
  & = \| [ \y^k + \gamma ( \tma_k x^k - \tmd ) ]_+
  \|^2 - \| \y^k \|^2 \nonumber\\
  & \leq \| \y^k + \gamma ( \tma_k x^k - \tmd ) \|^2 -
  \| \y^k \|^2 \nonumber\\
  & = - 2 \gamma \langle \tmd - \tma_k x^k, \y^k \rangle +
  \gamma^2 \| \tma_k x^k - \tmd \|^2 \nonumber\\
  & \leq - 2 \gamma \langle \tmd - \tma_k x^k, \y^k \rangle + m
  ( \bar{a} + \du )^2 \gamma^2, \nonumber
\end{align}
where the last inequality again uses the relation $\gamma^2 \| \tma_k x^k
- \tmd \|^2 \leq m ( \bar{a} + \du )^2 \gamma^2$. Then we take expectation and telescope over $k = 1, \ldots, n$ to obtain
\begin{align}
  \sum_{k = 1}^n \mathbb{E} [ \|\y^{k + 1}\|^2 - 
  \| \y^k \|^2 ] & ={}\mathbb{E} [ \|\y^{n + 1}\|^2
  ] - \|\y^1\|^2\nonumber\\
  & \leq{} \sum_{k = 1}^n \mathbb{E} [ - 2 \gamma \langle \tmd -
  \tma_k x^k, \y^k \rangle + m ( \bar{a} + \du )^2 \gamma^2
  ] . \nonumber
\end{align}
Re-arranging the terms, we have
\begin{align}
  \mathbb{E} \Big[\sum_{k = 1}^n 2 \gamma \langle \tmd - \tma_k x^k, \y^k \rangle \Big]
  \leq \sum_{k = 1}^n m ( \bar{a} + \du )^2 \gamma^2 -\mathbb{E}
  [ \| \y^{n + 1} \|^2 ]
  \leq{} &  mn ( \bar{a} + \du )^2 \gamma^2, \nonumber
\end{align}
Last we observe that \cite{10.48550/arxiv.2003.02513}
\begin{align}
  \sum_{k = 1}^n \mathbb{E} \left[ \frac{\hat{F}_{n - k + 1}^{\ast}}{n - k
  + 1} - c_k x^k \right]
  \leq{} \sum_{k = 1}^n \mathbb{E} [ \langle \tmd - \tma_k
  x^k,
  \y^k \rangle ] 
  \leq{} & \frac{m ( \bar{a} + \du )^2 \gamma n}{2} \nonumber
\end{align}
and plugging the bound back completes the proof.
\subsubsection{Proof of Theorem \ref{thm:explicit-theorem}}
We have, for some $\Delta$ independent of $\gamma$, that
\begin{align*}
 \mathbb{E} [ \rho ( \hat{\x} ) + v ( \hat{\x} ) ] & \leq 
  \Delta + \frac{m ( \bar{a} + \du )^2 \gamma n}{2} +
  \frac{\bar{c}}{\gamma \dl}\\
(\text{Taking }\gamma = \sqrt{\tfrac{2 \bar{c}}{\dl (
\bar{a} + \du )^2 m n}}) & = \Delta + 2 \left( \frac{( \bar{a} + \du )^2 \bar{c}}{2
  \dl} \right)^{1 / 2} \sqrt{m n},
\end{align*}
where $\frac{m ( \bar{a} + \du )^2
\gamma n}{2} = \frac{\bar{c}}{\gamma \dl}$ minimizes the right hand side and this completes the proof.

\subsection{Proof of Results in Section \ref{sec:theory-implicit}}
\label{app:proof-3-2}
In this section, we prove the results for implicit update. First we recall the update of \textbf{Algorithm \ref{alg:online}} under the implicit update.

\textbf{Update of Online Implicit Update}
\begin{align}
  \y^{k + 1} & = \argmin_{\y \geq \0}  \Big\{ \langle \tmd, \y \rangle+ [ c_k -
  \langle \tma_k, \y \rangle ]_+ + \frac{1}{2 \gamma} \|
  \y - \y^k \|^2 \Big\} \nonumber\\
  & = \argmin_{\y \geq \0} \Big\{ \langle \tmd, \y \rangle + s + \frac{1}{2 \gamma}
  \| \y - \y^k \|^2 : s \geq c_k - \langle \tma_k, \y
  \rangle, s \geq 0 \Big\} \nonumber\\
  x^k & = \lambda ( s \geq c_k - \langle \tma_k, \y \rangle
  ), \nonumber
\end{align}
where $\lambda ( s \geq c_k - \langle \tma_k, \y \rangle
)$ denotes the Lagrangian multiplier of the constraint $s \geq c_k -
\langle \tma_k, \y \rangle$.  As in the explicit case, we track the dual solution. 
\begin{lem}[Tracking the dual iteration]
  {\tmem{ \label{lem:implicit-dual-iter}
  
Under \ref{A1} and \ref{A2}, letting $\{ \y^k \}$ be the sequence of dual iterates generated by
  \textbf{Algorithm \ref{alg:online}} with implicit update and $\gamma_k \equiv \gamma$, then
  \begin{align}
    \| \y^{k + 1} - \y^k \| & \leq  \sqrt{m} ( \bar{a} +
    \du ) \gamma \nonumber\\
    \| \y^k \| & \leq \frac{m ( \bar{a} + \du )^2
    \gamma}{\dl} + \sqrt{m} ( \bar{a} + \du 
    ) \gamma  + \frac{\bar{c}}{\dl} . \nonumber
  \end{align}}}
\end{lem}
After bounding the dual iterations, we can move on to tracking constraint
violation and establish a bound of optimality gap, thus giving \textbf{Lemma \ref{lem:implicit-bound}} and \textbf{Theorem \ref{thm:implicit-theorem}}.
\subsubsection{Proof of Lemma \ref{lem:implicit-dual-iter}}
First we prove that $\| \y^{k + 1} - \y^k \| \leq \gamma \sqrt{m}
( \bar{a} + \du )$. Note that 
$$g ( \y ) = \langle \tmd,
\y \rangle + [ c_k - \langle \tma_k, \y \rangle
]_+ + \delta_{\y \geq \0}$$ is convex, where $\delta_{\y \geq \0}$ denotes the indicator function of $\mathbb{R}_+^m$. Then we invoke three-point lemma and
\[ \langle \tmd, \y^{k + 1} \rangle + [ c_k - \langle
   \tma_k, \y^{k + 1} \rangle ]_+ + \frac{1}{2 \gamma} \|
   \y^{k + 1} - \y^k \|^2 \leq \langle \tmd, \y^k \rangle +
   [ c_k - \langle \tma_k, \y^k \rangle ]_+ - \frac{1}{2
   \gamma} \| \y^{k + 1} - \y^k \|^2 . \]
Re-arranging the terms, we successively deduce that
\begin{align}
  \gamma^{- 1} \| \y^{k + 1} - \y^k \|^2 & \leq \langle \tmd,
  \y^k - \y^{k + 1} \rangle + [ c_k - \langle \tma_k, \y^k
  \rangle ]_+ - [ c_k - \langle \tma_k, \y^{k + 1}
  \rangle ]_+ \nonumber\\
  & \leq \| \tmd \| \cdot \| \y^k - \y^{k + 1} \| +
  \left| \langle \tma_k, \y^k - \y^{k + 1} \rangle \right|
  \label{eqn:implicit-1}\\
  & \leq \sqrt{m} \du \| \y^k - \y^{k + 1} \| + \sqrt{m} \bar{a}
  \| \y^k - \y^{k + 1} \| \label{eqn:implicit-2} \\
  & = \sqrt{m} ( \bar{a} + \du ) \| \y^k - \y^{k + 1}
  \|, \nonumber
\end{align}
where \eqref{eqn:implicit-1} uses Cauchy's inequality $\langle \tma, \tmb
\rangle \leq \| \tma \| \cdot \| \tmb \|$ and the relation
$[x]_+ - [y]_+ \leq | x - y |_+$; \eqref{eqn:implicit-2} again applies
Cauchy's inequality together with \ref{A1}, \ref{A2}. Dividing both
sides of the inequality by $\| \y^k - \y^{k + 1} \|$ shows $\|
\y^{k + 1} - \y^k \| \leq \gamma \sqrt{m} ( \bar{a} + \du )$.\condsep

Next we bound $\| \y^k \|$. Due to the complication of implicit
update, we resort to a constrained smooth formulation of the proximal
subproblem
\begin{eqnarray*}
  \min_{\y} & \langle \tmd, \y \rangle + s + \frac{1}{2 \gamma} \| \y - \y^k
  \|^2 & \text{Dual}\\
  \text{subject to} & s \geq c_k - \langle \tma_k, \y \rangle & x\\
  & s \geq 0, \y \geq \0 & v, \w,
\end{eqnarray*}
which is a convex quadratic programming problem. Now we check the Lagrangian function
\[ L ( \y, x, \s, x, v ) = \langle \tmd, \y \rangle + s
   + \frac{1}{2 \gamma} \| \y - \y^k \|^2 + x ( c_k -
   \langle \tma_k, \y \rangle - s ) - v s - \langle \y,
   \w \rangle, \]
where $\w$ is the multiplier of $\y$ and $v$ is the multiplier of $s$.
Writing the KKT conditions, we have
\begin{align}
  s & \geq c_k - \langle \tma_k, \y^{k + 1} \rangle \nonumber \\
  s & \geq 0 \nonumber\\
  \tmd + \gamma^{- 1} ( \y^{k + 1} - \y^k ) - \tma_k x^k - \w & =
  \0 \label{KKT-3}\\
  v + x^k & = 1 \nonumber\\
  \langle \y^{k + 1}, \w \rangle & = 0\nonumber\\
  x^k ( c_k - \langle \tma_k, \y^{k + 1} \rangle - s ) &
  = 0 \label{KKT-6} \\
  v s & = 0 \nonumber\\
  ( x^k, v, \y^{k + 1}, \w ) & \geq \0 \nonumber
\end{align}

and from \eqref{KKT-3} we know that
\begin{align}
	\y^{k + 1} = \y^k - \gamma ( \tmd - \tma_k x^k ) + \gamma \w \label{eqn:implicit-update-formula}
\end{align}
for some $\w \geq \0$. Also we notice that since $\langle \y^{k + 1}, \w
\rangle = 0$, $w_i = 0$ whenever $y_i^{k + 1} > 0$, which implies
\begin{align}
  \| \y^{k + 1} \|^2 & = \| \y^k - \gamma ( \tmd - \tma_k
  x^k ) + \gamma \w \|^2 \nonumber\\
  & \leq \| \y^k - \gamma ( \tmd - \tma_k x^k ) \|^2
  \nonumber\\
    & = \| \y^k \|^2 - 2 \gamma \langle \tmd - \tma_k x^k, \y^k
  \rangle + \gamma^2 \| \tmd - \tma_k x^k \|^2 \label{eqn:lem-6-2} \\
  & \leq \| \y^k \|^2 + 3 m ( \bar{a} + \du )^2
  \gamma^2 + 2 \gamma \bar{c} - 2 \gamma \dl \| \y^k \| \label{eqn:implicit-descent}. 
\end{align}
where \eqref{eqn:implicit-descent} is again by Cauchy's inequality and \ref{A1}, \ref{A2}, $- 2
\gamma \langle \tmd, \y^k \rangle \leq 2 \gamma \dl \| \y^k
\|_1 \leq - 2 \gamma \dl \| \y^k \|$, and that
\begin{align}
  \langle \tma_k x^k, \y^k \rangle & = \langle \tma_k x^k,
  \y^{k + 1} \rangle - \langle \tma_k x^k, \y^{k + 1} - \y^k
  \rangle \nonumber\\
  & \leq \langle \tma_k x^k, \y^{k + 1} \rangle + \| \tma_k
  \| \cdot \| \y^{k + 1} - \y^k \| \cdot |x^k| \nonumber\\
  & \leq \bar{c} + \gamma m \bar{a} ( \bar{a} + \du ) \label{eqn:lem-6-1} \\
  & \leq \bar{c} + \gamma m ( \bar{a} + \du )^2, \nonumber
\end{align}
where \eqref{eqn:lem-6-1} uses $\langle \tma_k x^k, \y^{k + 1} \rangle \leq \bar{c}$ from \eqref{KKT-6} and we invoke the bound $\| \y^{k + 1} - \y^k \| \leq \gamma \sqrt{m} ( \bar{a} + \du )$.

On the other hand, we know that
\begin{align}
  \| \y^{k + 1} \| & = \| \y^{k + 1} - \y^k + \y^k \|
  \nonumber\\
  & \leq \| \y^{k + 1} - \y^k \| + \| \y^k \|
  \nonumber\\
  & \leq \| \y^k \| + \gamma \sqrt{m} ( \bar{a} + \du )\label{eqn:implicit-4},
\end{align}
where \eqref{eqn:implicit-4} again uses $\|
\y^{k + 1} - \y^k \| \leq \gamma \sqrt{m} ( \bar{a} + \du )$.
By exactly the same argument as in \tmtextbf{Lemma \ref{lem:explicit-dual-iter}}, we know that
$$\| \y^k \| \leq \frac{3m ( \bar{a} + \du )^2
\gamma}{\dl} +  \sqrt{m} ( \bar{a} + \du
) \gamma + \frac{\bar{c}}{\dl}$$ and this completes the proof.

\subsubsection{Proof of Lemma \ref{lem:implicit-bound}}
First we consider constraint violation and recall that we run implicit update \eqref{eqn:implicit-update-formula}
\[ \y^{k + 1} = \y^k - \gamma ( \tmd - \tma_k x^k ) + \gamma \w
   \geq \y^k - \gamma ( \tmd - \tma_k x^k ). \]
Hence telescoping gives
\begin{align}
  \mathbb{E} [ \| \A \hat{\x} - \tmb \|_+ ]  &
  \leq  \gamma^{- 1} \mathbb{E} [ \| \y^{k + 1} \| ]
  \nonumber
  \leq \frac{3m ( \bar{a} + \du )^2}{\dl} + \sqrt{m} (
  \bar{a} + \du ) + \frac{\bar{c}}{\gamma \dl} . \nonumber
\end{align}
As for the optimality gap, we have, similar to \cite{10.48550/arxiv.2003.02513}, that
\[ \sum_{k = 1}^n \mathbb{E} \Big[ \frac{\hat{F}_{n - k + 1}^{\ast}}{n - k +
   1} - c_k x^k \Big] \leq \sum_{k = 1}^n \mathbb{E} [ \langle
   \tmd, \y^k \rangle + [ c_k - \langle \tma_k, \y^k
   \rangle ]_+ - c_k x^k ] \]
Now we look again into the KKT conditions, which, after simplification, gives
\begin{align}
  s & \geq [ c_k - \langle \tma_k, \y^{k + 1} \rangle
  ]_+ \label{eqn:implicit-5}\\
  x^k ( c_k - \langle \tma_k, \y^{k + 1} \rangle - s ) &
  = 0 \label{eqn:implicit-6}\\
  (1 - x^k) s & = 0. \label{eqn:implicit-7}
\end{align}
Combining the above relations, we successively deduce that
\begin{align}
  c_k x^k & = ( \langle \tma_k, \y^{k + 1} \rangle + s )
  x^k \label{eqn:lem-2-1}\\
  & = \langle \tma_k, \y^{k + 1} \rangle x^k + s x^k\label{eqn:lem-2-3} \\
  & = \langle \tma_k, \y^{k + 1} \rangle x^k + s \label{eqn:lem-2-3}\\
  & \geq \langle \tma_k, \y^{k + 1} \rangle x^k + [ c_k -
  \langle\tma_k, \y^{k + 1} \rangle ]_+, \nonumber
\end{align}
where \eqref{eqn:lem-2-1} re-arranges \eqref{eqn:implicit-6} and \eqref{eqn:lem-2-3} uses \eqref{eqn:implicit-7}. Plugging it back, we can derive the following bound.
\begin{align}
  & \sum_{k = 1}^n \mathbb{E} [ \langle \tmd, \y^k \rangle +
  [ c_k - \langle \tma_k, \y^k \rangle ]_+ - c_k x^k
  ] \nonumber\\
  \leq{} & \sum_{k = 1}^n \mathbb{E} [ \langle \tmd, \y^k
  \rangle + [ c_k - \langle \tma_k, \y^k \rangle
  ]_+ - \langle \tma_k, \y^{k + 1} \rangle x^k - [ c_k -
  \langle \tma_k, \y^{k + 1} \rangle ] _+] \nonumber\\
  ={} & \sum_{k = 1}^n \mathbb{E} [ \langle \tmd - \tma_k x^k, \y^k
  \rangle ] + \sum_{k = 1}^n \mathbb{E}  [ \langle
  \tma_k x^k, \y^k - \y^{k + 1} \rangle ] + \sum_{k = 1}^n
  \mathbb{E} \{ [ c_k - \langle \tma_k, \y^k \rangle
  ]_+ - [ c_k - \langle \tma_k, \y^{k + 1} \rangle
  ]_+ \} \nonumber
\end{align}
Next we bound the last two summations by
\begin{align}
  \sum_{k = 1}^n \mathbb{E}  [ \langle \tma_k x^k, \y^k - \y^{k + 1}
  \rangle ] & \leq \sum_{k = 1}^n \mathbb{E}  [ \|
  \tma_k \| \cdot \| \y^k - \y^{k + 1} \| ]
  \leq m n \gamma \bar{a} ( \bar{a} + \du ) \nonumber
\end{align}
\begin{align}
\sum_{k = 1}^n \mathbb{E} \{ [ c_k - \langle \tma_k, \y^k
  \rangle ]_+ - [ c_k - \langle \tma_k, \y^{k + 1}
  \rangle ]_+ \}
   & \leq \sum_{k = 1}^n \mathbb{E} [ | \langle \tma_k, \y^k -
  \y^{k + 1} \rangle | ] \nonumber\\
  &\leq  \sum_{k = 1}^n \mathbb{E}  [ \| \tma_k \| \cdot
  \| \y^k - \y^{k + 1} \| ] \nonumber\\
  &  \leq m n \gamma \bar{a} ( \bar{a} + \du ) \nonumber
\end{align}
with Cauchy's inequality, \textbf{Lemma \ref{lem:implicit-dual-iter}} and \ref{A2}.
Then we re-arrange \eqref{eqn:lem-6-2} and obtain
\begin{align}
  & 2 \gamma \langle \tmd - \tma_k x^k, \y^k \rangle \nonumber\\
  \leq & ~ \| \y^k \|^2 - \| \y^{k + 1} \|^2 + \gamma^2
  \| \tmd - \tma_k x^k \|^2 \nonumber\\
  \leq & ~ \| \y^k \|^2 - \| \y^{k + 1} \|^2 + \gamma^2 m
  ( \bar{a} + \du )^2 \label{eqn:implicit-3}
\end{align}
and \eqref{eqn:implicit-3} uses \ref{A1} and \ref{A2} and the fact that
$|x^k| \leq 1$.  Finally, we telescope
\[ 2 \gamma \langle \tmd - \tma_k x^k, \y^k \rangle \leq \|
   \y^k \|^2 - \| \y^{k + 1} \|^2 + \gamma^2 m ( \bar{a}
   + \du )^2, \]
as in the previous analysis to get
\[ \sum_{k = 1}^n \mathbb{E} [ \langle \tmd - \tma_k x^k, \y^k
   \rangle ] \leq \frac{m ( \bar{a} + \du )^2 n   \gamma}{2}
 . \]
Putting all the bounds together, we have
\[ \Ebb[\rho ( \hat{\x} ) ]\leq \frac{m n ( \bar{a} + \du )^2}{2}
   \gamma + 2 m n \gamma \bar{a} ( \bar{a} + \du ) \leq \frac{5 m
   ( \bar{a} + \du )^2 n \gamma}{2} \]
and this completes the proof.

\subsubsection{Proof of Theorem \ref{thm:implicit-theorem}}
The proof works exactly in the same way as \tmtextbf{Theorem \ref{thm:explicit-theorem}} by
observing that
\[ \Ebb[\rho ( \hat{\x} ) + v ( \hat{\x} )] \leq \Delta + \frac{5 m n
   ( \bar{a} + \du )^2}{2} \gamma + \frac{\bar{c}}{\gamma \dl} \]
and taking optimal $\gamma^{\ast} = \sqrt{\frac{2 \bar{c}}{5 \dl (
\bar{a} + \du )^2 m n}}$ to minimize the right hand side.

\subsection{Proof of Results in Section \ref{sec:theory-duplicate}}
\label{app:proof-3-3}
\subsection{Proof of Theorem \ref{thm:dup}}
In this section we consider the variable duplication scheme. Given an LP
\begin{eqnarray*}
  \max_{\x} & \langle \tmc, \x \rangle & \\
  \text{subject to} & \A \x \leq \tmb & \\
  & \0 \leq \x \leq \1 & 
\end{eqnarray*}
and its duplicated version
\begin{eqnarray*}
  \max_{\left\{ \x_j \right\}} & \sum_{j = 1}^K
  \langle \tmc, \x_j \rangle & \\
  \text{subject to} & \sum_{j = 1}^K \A \x_j \leq K \tmb & \\
  & \0 \leq \x_j \leq \1, & 
\end{eqnarray*}
It's clear that the duplicated LP also satisfies \ref{A1},
\ref{A2}. Also we know that \ref{A3} can be satisfied by an arbitrarily small perturbation of the objective coefficients. Then we immediately have $\sum_{j  = 1}^K \langle \tmc, \x_j^{\ast} \rangle
= K \langle \tmc, \x^{\ast} \rangle$ up to some arbitrarily small perturbation.
Suppose that we apply \tmtextbf{Theorem \ref{thm:explicit-theorem}} or \tmtextbf{Theorem \ref{thm:implicit-theorem}} to get
$\{ \x_j' \}$ such that
\begin{align}
  \Ebb [K \langle \tmc,
  \x^{\ast} \rangle - \sum_{j  = 1}^K \langle \tmc, \x'_j \rangle] & =\mathcal{O} ( m \log n+ \sqrt{n K} \log n + \sqrt{m n
  K} + \sqrt{n K} \log K) \nonumber\\
  \Ebb \Big[\Big\| \Big[ \sum_{j  = 1}^K \A \x_j' - K \tmb \Big]_+ \Big\|\Big] &
  =\mathcal{O} ( m + \sqrt{m n K} ) . \nonumber
\end{align}
Then in view of $\hat{\x} = \frac{1}{K} \sum_{j  = 1}^K \x_j'$ we have
\[ \Ebb[\rho ( \hat{\x} )] = \Ebb\Big[ \langle \tmc, \x^{\ast} \rangle - \langle \tmc, \frac{1}{K} \sum_{k
   = 1}^K \x'_k \rangle \Big]
   =\mathcal{O} \left( \frac{m \log n}{K}+ \sqrt{\frac{n}{K}} \log n + \sqrt{\frac{m n}{K}} + \sqrt{\frac{n}{K}} \log K
   \right) \]
\[ \Ebb[v ( \hat{\x} )] = \Ebb\Big[ \Big\| \Big[ \A \Big( \frac{1}{K}
   \sum_{j  = 1}^K \x_j' \Big) - \tmb \Big]_+ \Big\|\Big] =\mathcal{O} \left(
   \frac{m}{K}+\sqrt{\frac{m n}{K}} \right) \]
assuming that $K = \mathcal{O}(n)$, and this completes the proof.

\subsection{Additional Experiment on Violated Assumption} \label{app:addtional-exp}

In this section, we carry out additional experiment to see what happens when some of our assumptions
are (nearly) violated. Specifically we consider the MKP problem with $\mathbf{b} = \textbf{1}$.
\begin{eqnarray*}
  \max_{\mathbf{x}} & \langle \mathbf{c}, \mathbf{x} \rangle & \\
  \text{subject to} & \mathbf{Ax} \leq \textbf{1} & \\
  & \textbf{0} \leq \mathbf{x} \leq \textbf{1}, & 
\end{eqnarray*}
Note that when $n$ is large, \tmtextbf{A1} becomes asymptotically violated since $\underline{d} = \frac{1}{n}$ becomes closer to 0.\\

\paragraph{Testing Configuration and Setup} We configure the algorithm as follows (\textbf{Section \ref{sec:experiments}} gives a more detailed description).
\begin{enumerate}[label=\textbf{\arabic*).},ref=\rm{\textbf{A\arabic*}},leftmargin=*,]
  \item \tmtextbf{Dataset}. We take $(m, n) \in \{ (5, 100), (8, 1000), (16,
  2000), (32, 4000) \}, \sigma = 1$.
  
  \item \tmtextbf{Initial Point}. We let online algorithms start from
  {\tmstrong{0}}.
  
  \item \tmtextbf{Feasibility}. We force the algorithm to respect constraint
  violation.
  
  \item \tmtextbf{Duplication}. We allow $K \in \{ 1, 2, 4, 8, 16, 32, 64,
  128 \}$
  
  \item \tmtextbf{Stepsize}. We take $\gamma = (K m n)^{- 1 / 2}$.
\end{enumerate}
\begin{figure}[h]
\centering
{\includegraphics[scale=0.22]{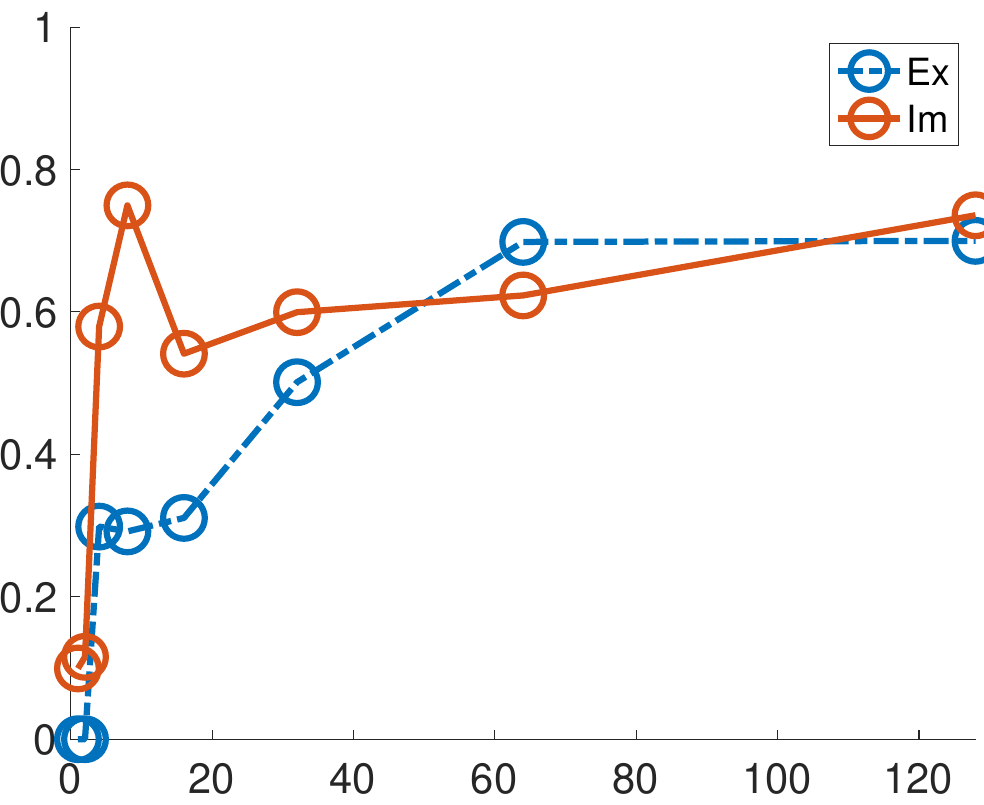}}{\includegraphics[scale=0.22]{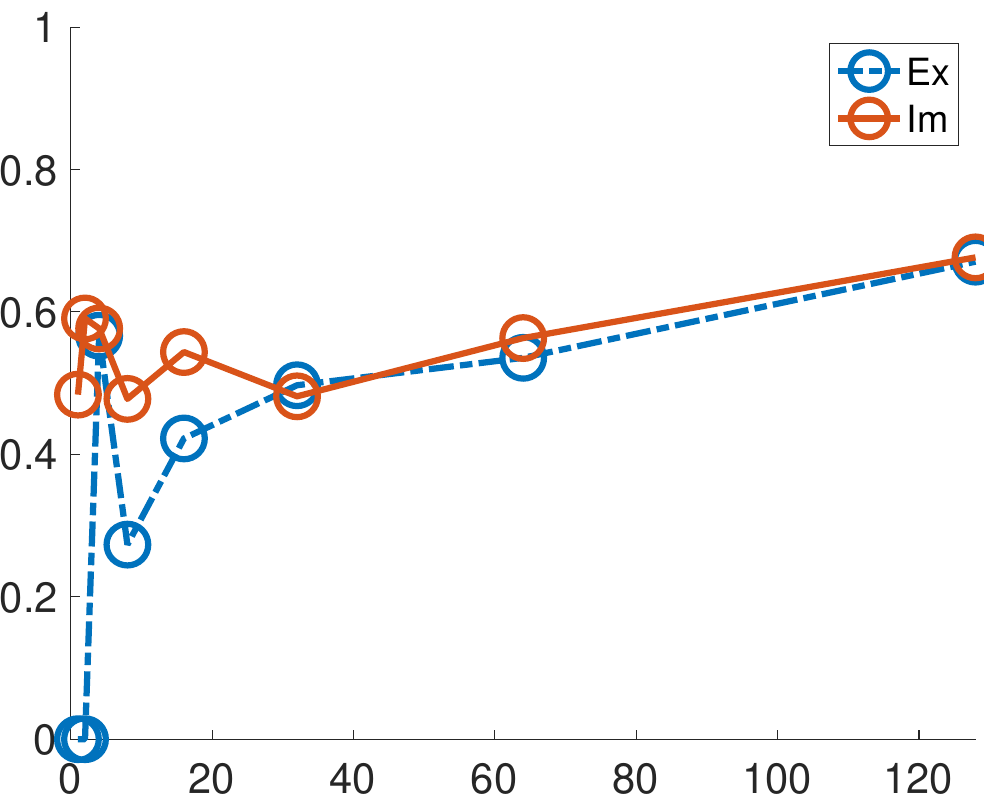}}{\includegraphics[scale=0.22]{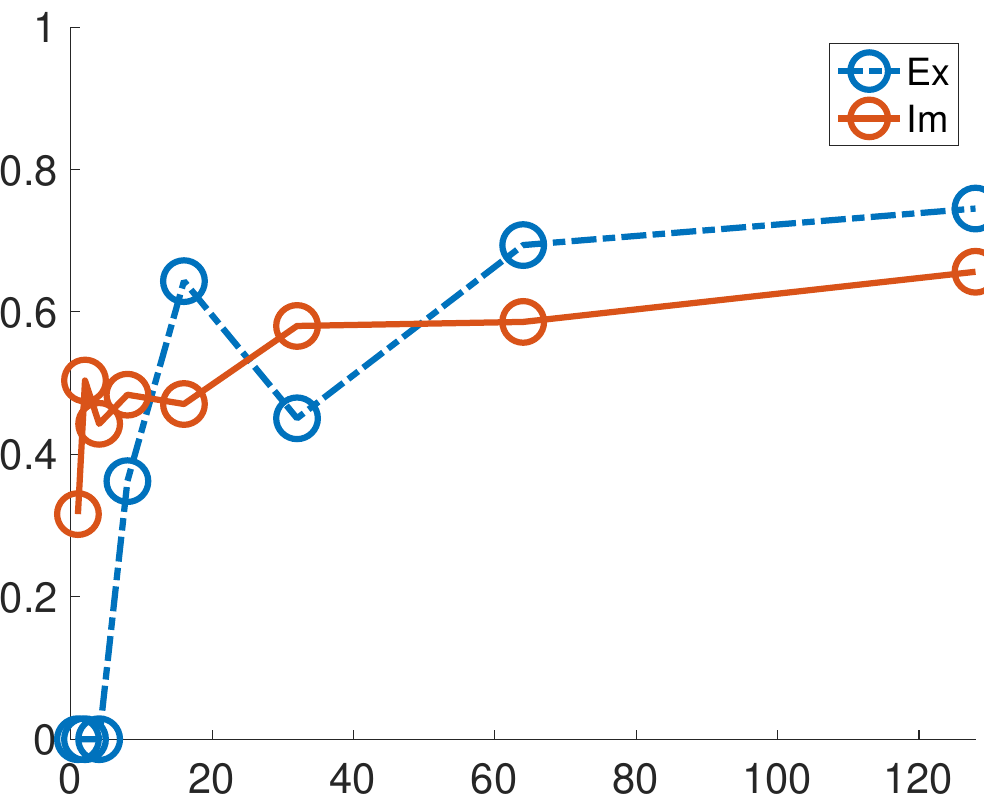}}{\includegraphics[scale=0.22]{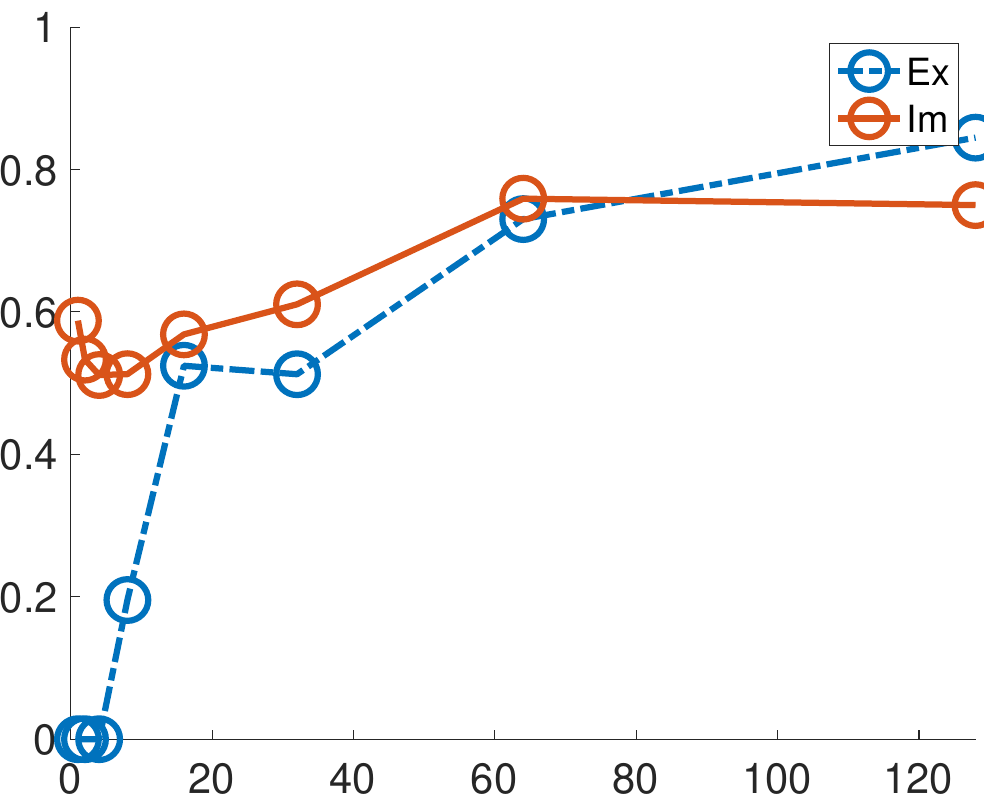}}
\caption{From left to right: $(m, n) \in \{ (5, 100), (8, 1000), (16, 2000),
  (32, 4000) \}$
  The axis represents $K$ parameter ranging from 1 to 128. y-axis, relative
  optimality gap}
\end{figure}

We see that as suggested by our theory, as $\underline{d} = \frac{1}{n}$ gets
smaller, online algorithms perform poorly when $K = 1$, but when we increase
$K$, both explicit and implicit update gradually retrieve better solutions.

\newpage
\section{Practical Aspects of Online Linear Programming}\label{app:applications}
In this section, we discuss the practical aspects of online linear programming in more details. In \textbf{Section \ref{app:impl}} and \textbf{Section \ref{app:subproblem}}, we show when and how online updates can be carried out efficiently; In \textbf{Section \ref{app:subproblem}}, we apply online explicit update to more LP problem types and compare its performance with other LP solvers; In \textbf{Section \ref{app:sifting}}, we formalize the contents of LP sifting and present more implementation details. Finally we show in \textbf{Section \ref{app:mip}} that online LP is applicable to speedup mixed integer programming.
\subsection{Fast $\mathcal{O}(\text{nnz}(\A))$ Implementation For the Explicit Update} \label{app:impl}
In this section, we discuss the practical aspects of implementing the explicit update.
The computation of explicit update comes from the following two steps
\begin{align}
  x^k \leftarrow{} & \mathbb{I} \{ c_k > \langle \tma_k, \y
  \rangle \} \nonumber\\
  \y^{k + 1} \leftarrow{} & [ \y^k + \tma_k x^k - \tmd ]_+, \nonumber
\end{align}
where we take $\gamma = 1, d_i \equiv d$ without loss of generality. Due to the
\texttt{axpy} operation $\y + \tma_k x^k - \tmd$ and
the full density of $\tmd$, a direct implementation of \texttt{axpy} and
$[ \x ]_+$ will immediately result in $\mathcal{O} (m n)$ flops and
is undesirable. To address this, we define $p_j (k) \assign \max_l \{ l < k :
a_{l, j} \neq 0 \}$, the latest iteration $l$ where $a_{l,j}$ is nonzero. Then we
observe that
\begin{align}
  y^{k + 1}_j ={} & [y^k_j - a_{k, j} x^k - d]_+ \nonumber\\
  ={} & [[y^{k - 1}_j - a_{k - 1, j} x^{k - 1} - d]_+ - a_{k, j} x^k - d]_+
  \nonumber\\
  ={} & \cdots \nonumber\\
  ={} & [[y_j^{p (k)} - a_{p (k), j} x^{p_j (k)} - d]_+ - (k - p (k)) d]_+,
  \nonumber
\end{align}
where we use the relation $[[a - b]_+ - b]_+ = [a - 2 b]_+, b \geq 0$ recursively.
Further we observe that
\[ a_j^{k + 1} y^{k + 1}_j = a_j^{k + 1} [[y_j^{p (k)} - a_{p (k), j} x^{p_j
   (k)} - d]_+ - (k - p (k)) d]_+, \]
which implies that we only need to evaluate $[[y_j^{p (k)} - a_{p (k), j}
x^{p_j (k)} - d]_+ - (k - p (k)) d]_+$ if $a_j^{k + 1} \neq 0$ and the
operation takes $\mathcal{O} (1)$. Hence we have overall complexity of
$\mathcal{O}( \text{nnz}( \A ))$ to compute all the
primal estimates and $\mathcal{O} (m)$ to recover $\y^{n + 1}$. In
practice we can maintain an $\mathbb{N}^m$ array \texttt{lastUpdate[m]} to
record and to update $p_j (k)$ for each $j \in [m]$.
\subsection{Cases of Simple Implicit Update} \label{app:subproblem}

In this section we focus on the practical implementation of the online implicit update and discuss the case where the implicit is easier to compute. Here we consider
\begin{eqnarray*}
  \min_{\y} & \langle \tmd, \y \rangle + [ c - \langle
  \tma, \y \rangle ]_+ + \frac{1}{2\gamma} \| \y - \mathbf{z}
  \|^2 & \\
  \text{subject to} & \y \geq \0 . & 
\end{eqnarray*}

where $\tma \geq \0$. Let $\y^+$ be the optimal solution to the problem and we do case analysis.

\tmtextbf{Case 1}. $c - \langle \tma, \y^+ \rangle > 0$. Then $\y^+
= [\mathbf{z} - \gamma^{- 1} ( \tma - \tmd )]_+$.

{\tmstrong{Case 2}}. $c - \langle \tma, \y^+ \rangle < 0$. Then
$\y^+ = [\mathbf{z} - \gamma^{- 1} \tmd]_+$.

\tmtextbf{Case 3}. $c - \langle \tma, \y^+ \rangle = 0$. Then we
have, equivalently, that $\y^+$ is the optimal solution to the following
problem.
\begin{eqnarray*}
  \min_{\y} & \langle \tmd, \y \rangle + \frac{1}{2\gamma} \|
  \y - \mathbf{z} \|^2 & \\
  \text{subject to} & \langle \tma, \y \rangle = c & \\
  & \y \geq \0 . & 
\end{eqnarray*}
Since $\tma \geq \0$, letting $(A, \bar{A})$ be the partition of $[m]$ such
that $\tma_A > \0$ and we can separate the problem into
\begin{eqnarray*}
  \min_{\y_{\bar{A}}} & \langle \tmd_{\bar{A}}, \y_{\bar{A}}
  \rangle + \frac{1}{2\gamma} \| \y_{\bar{A}} - \mathbf{z}_{\bar{A}}
  \|^2 & \\
  \text{subject to} & \y_{\bar{A}} \geq \0 & 
\end{eqnarray*}
and
\begin{eqnarray*}
  \min_{\y_A} & \langle \tmd_A, \y_A \rangle + \frac{1}{2\gamma}
  \| \y_A - \mathbf{z}_A \|^2 & \\
  \text{subject to} & \langle \tma_A, \y_A \rangle = c & \\
  & \y_A \geq \0, & 
\end{eqnarray*}
where $\y_{\bar{A}}$ can be efficiently updated and $\y_A$ can be written as
\[ \y_A = \text{Proj}_{\Delta_{\tma_A / c}} ( \mathbf{z}_A - \gamma^{- 1}
   \tmd_A ), \]
where $\text{Proj}_{\Delta_{\tma_A}}$ denotes orthogonal projection onto
the weighted simplex
\[ \Delta_{c^{-1}\tma_A} \assign \{ \y_A : \langle c^{- 1} \tma_A, \y_A
   \rangle = 1, \y_A \geq \0 \} \]
and can be efficiently solved using sorting-based method proposed in
{\cite{perez2022efficient}}.

\begin{rem}
	In practice, we can try the solutions from the first two cases and verify if they satisfy the conditions from case analysis.  
	If in both cases the condition is violated, we invoke the above subroutine and compute the projection.
\end{rem}

\subsection{Direct Approximate LP Solving} \label{app:direct}

In this section, we evaluate the performance of online LP solver for direct LP solving. The experiments in this section consist of two parts. First we focus
on the CPU running time of our proposed method, and then we turn to more LP
instances as an extension of the experiments in {\tmstrong{Section \ref{sec:experiments}}}.

\paragraph{CPU Time Evaluation}

First we evaluate the CPU running time of our methods. As we already
discussed, we can implement the online LP algorithm in $\mathcal{O}
(\text{nnz} (\mathbf{A}))$ time, which implies our method is highly scalable.
The following table summarizes the total CPU time of {\textbf{Algorithm \ref{alg:online-dup}}} with $K=100$ under different settings of $m, n$ and $\text{nnz}(\mathbf{A})$. Now that explicit update is independent of specific numerical values of the matrices, we use the same way of data generation as in {\tmstrong{Section \ref{sec:experiments}}}.\condsep

\begin{table}[h]
\centering
  \caption{CPU Time evaluation of explicit update. \label{table:cputime} Time given in CPU seconds.}
  \begin{tabular}{cccccccccccc}\toprule
    $m$ & $n$ & nnz & Time & $m$ & $n$ & nnz & Time & $m$ & $n$ & nnz & Time\\
     \midrule
    $10^2$ & $10^2$ & $10^3$ & 0.00 & $10^3$ & $10^2$ & $10^4$ & 0.00 & $10^4$
    & $10^2$ & $10^4$ & 0.01\\
    $10^2$ & $10^3$ & $10^4$ & 0.00 & $10^3$ & $10^3$ & $10^4$ & 0.00 & $10^4$
    & $10^3$ & $10^5$ & 0.04\\
    $10^2$ & $10^4$ & $10^4$ & 0.02 & $10^3$ & $10^4$ & $10^5$ & 0.05 & $10^4$
    & $10^4$ & $10^6$ & 0.38\\
    $10^2$ & $10^5$ & $10^5$ & 0.26 & $10^3$ & $10^5$ & $10^6$ & 0.51 & $10^4$
    & $10^5$ & $10^5$ & 0.28\\
    $10^2$ & $10^6$ & $10^6$ & 2.61 & $10^3$ & $10^6$ & $10^5$ & 0.72 & $10^4$
    & $10^6$ & $10^6$ & 2.77\\
    \bottomrule
    &  &  &  &  &  &  &  &  &  &  & \\
    \toprule
    $m$ & $n$ & nnz & Time & $m$ & $n$ & nnz & Time &  &  &  & \\
     \midrule
    $10^5$ & $10^2$ & $10^5$ & 0.09 & $10^6$ & $10^2$ & $10^6$ & 1.54 & 
    &  &  & \\
    $10^5$ & $10^3$ & $10^6$ & 0.57 & $10^6$ & $10^3$ & $10^5$ & 0.57 & 
    &  &  & \\
    $10^5$ & $10^4$ & $10^5$ & 0.11 & $10^6$ & $10^4$ & $10^6$ & 1.61 & 
    &  &  & \\
    $10^5$ & $10^5$ & $10^6$ & 0.79 & $10^6$ & $10^5$ & $10^7$ & 11.01 &
    &  &  & \\
    $10^5$ & $10^6$ & $10^7$ & 8.08 & $10^6$ & $10^6$ & $10^8$ & 120.00 &
    &  &  &\\
    \bottomrule 
  \end{tabular}
\end{table}

\textbf{Table \ref{table:cputime}} suggests, for the same $n$, CPU time our method increases almost linearly with respect to $\text{nnz} (\mathbf{A})$, which verifies that the implementation from \textbf{Section \ref{app:impl}} is highly scalable to huge linear programs.

\paragraph{Real and synthetic LPs}

Now we switch to a more practical setting where we employ online algorithms to
solve real-life instances and compare performance of online algorithms with LP
solvers. Our setup is given as follows.\condsep

\tmtextbf{Testing configuration and setup}
\begin{enumerate}[label=\textbf{\arabic*).},ref=\rm{\textbf{A\arabic*}},leftmargin=*,]
  \item \tmtextbf{Dataset}. Our dataset comes from three sources.
  \tmtextbf{1)}. 7 instances that are LP relaxations from \texttt{MIPLIB}.
  \tmtextbf{2)}. 4 MKP instances generated according to the statistics in
  \tmtextbf{Table \ref{table:cputime}}. \tmtextbf{3)}. Modified \texttt{Netlib}
  instances.
  
  \item \tmtextbf{Initial point}. We let online algorithm start from
  {\tmstrong{0}}
  
  \item \tmtextbf{Feasibility}. We do not enforce feasibility of the constraints
  
  \item \tmtextbf{Duplication}. We allow a maximum of $K = 5000$ variable
  duplications
  
  \item \tmtextbf{Stopping criterion}. We let the algorithm stop if
  \[ \max \left\{ \frac{\| [\mathbf{Ax} - \mathbf{b}]_+ \|}{\| \mathbf{b} \|_1
     + 1}, \frac{\langle \mathbf{b}^{\top} \mathbf{y} \rangle + \langle
     \mathbf{u}, [\mathbf{c} - \mathbf{A}^{\top} \mathbf{y}]_+\rangle - \left\langle
     \mathbf{c},{{\tmstrong{x}}} \right\rangle}{| \langle \mathbf{b}^{\top}
     \mathbf{y} \rangle + \langle \mathbf{u}, [\mathbf{c} - \mathbf{A}^{\top}
     \mathbf{y}]_+) | + \left| \left\langle \mathbf{c},{{\tmstrong{x}}}
     \right\rangle \right| + 1} \right\} \leq \varepsilon = 5 \times 10^{- 3}
  \]
  \item {\tmstrong{Modification of \texttt{Netlib} instances}}. Some
  \texttt{Netlib} instances do not meet assumptions for the online
  algorithm and thus prohibits direct solving. We modify
  \texttt{Netlib} instances by \tmtextbf{1)}. Taking $\hat{b}_i = \max \{
  b_i, 10^{- 3} \}$. \tmtextbf{2)}. Enforcing upperbound $u_i = \min \{ u_i,
  10^2 \}$. \tmtextbf{3}). Changing constraint senses into $\leq$.
\end{enumerate}

\begin{table}[h]
\centering
  \caption{Time of solving \texttt{MIPLIB} instances to an \texttt{1e-03} relative
  accuracy solution and comparison with  \texttt{Gurobi v9.5}. Synthetic MKP instance
  \texttt{mkp-\{i\}-\{j\}} stands for MKP instance of $10^i$ rows and $10^j$
  columns. GTime:  \texttt{Gurobi} solution time.}
  \begin{tabular}{ccccc|ccccc}\toprule
    Instance & pInf & Gap & Time & GTime & Instance & pInf & Gap & Time &
    GTime\\ \midrule
    \texttt{2club200v15p5scn} & 1.7e-04 & 5.5e-02 & 3.44 & 0.32 &
    \texttt{mkp-2-5} & 4.4e-03 & 4.2e-03 & 0.03 & 0.44\\
    \texttt{cdc7-4-3-2} & 5.6e-05 & 5.0e-03 & 1.27 & 24.00 &
    \texttt{mkp-2-6} & 5.3e-05 & 4.9e-03 & 0.17 & 0.22\\
    \texttt{cod105} & 1.0e-06 & 9.6e-03 & 1.29 & 0.99 &
    \texttt{mkp-2-7} & 5.3e-05 & 5.0e-03 & 0.62 & 5.70\\
    \texttt{p6b} & 1.0e-04 & 1.2e-05 & 0.09 & 0.16 & \texttt{mkp-3-5}
    & 1.8e-03 & 4.9e-03 & 0.06 & 0.11\\
    \texttt{m100n500k4r1} & 9.7e-05 & 3.4e-03 & 0.12 & 0.08 &
    \texttt{mkp-3-6} & 1.8e-03 & 4.9e-03 & 0.29 & 0.54\\
    \texttt{manna81} & 1.0e-04 & 4.9e-03 & 0.62 & 0.08 &
    \texttt{mkp-3-7} & 4.4e-03 & 2.7e-03 & 1.38 & 44.70\\
    \texttt{queens-30} & 4.0e-05 & 4.8e-03 & 0.22 & 0.63 &  &  &  &  & \\
    \bottomrule
  \end{tabular}
\end{table}

\begin{table}[h]
\centering
\caption{Experiment on \textit{modified} \texttt{Netlib} LP datasets. GTime:  \texttt{Gurobi} solution time}
\resizebox{0.7\textwidth}{!}{
\begin{tabular}{ccccc|ccccc}
  \toprule
    Instance & pInf & Gap & Time & GTime & Instance & pInf & Gap & Time &
    GTime\\
    \midrule
    25fv47 & 9.9e-05 & 4.94-03 & 0.45 & 0.02 & osa-30 & 0.0e+00 & 0.00+00 &
    0.00 & 0.13\\
    80bau3b & 6.4e-04 & 2.18-02 & 1.13 & 0.10 & osa-60 & 0.0e+00 & 0.00+00 &
    0.00 & 0.15\\
    adlittle & 2.8e-03 & 1.92-02 & 0.60 & 0.08 & pds-02 & 0.0e+00 & 2.41-05 &
    0.00 & 0.07\\
    afiro & 0.0e+00 & 4.92-03 & 0.01 & 0.07 & pds-06 & 4.1e-05 & 1.10-05 &
    0.00 & 0.13\\
    agg & 0.0e+00 & 0.00+00 & 0.00 & 0.07 & pds-10 & 2.4e-05 & 5.88-05 & 0.00
    & 0.14\\
    agg2 & 0.0e+00 & 0.00+00 & 0.00 & 0.07 & pds-20 & 1.7e-05 & 1.44-04 & 0.00
    & 0.12\\
    agg3 & 0.0e+00 & 0.00+00 & 0.00 & 0.07 & perold & 0.0e+00 & 1.49-02 & 0.69
    & 0.08\\
    bandm & 1.6e-03 & 4.57-02 & 0.62 & 0.08 & pilot.ja & 9.2e-05 & 4.25-04 &
    0.01 & 0.08\\
    beaconfd & 0.0e+00 & 0.00+00 & 0.00 & 0.07 & pilot & 4.9e-05 & 1.32-02 &
    1.14 & 0.13\\
    blend & 2.8e-04 & 1.63-03 & 0.60 & 0.08 & pilot.we & 0.0e+00 & 3.56-05 &
    0.00 & 0.07\\
    bnl1 & 8.2e-04 & 1.28-01 & 0.67 & 0.08 & pilot4 & 0.0e+00 & 1.07-02 & 0.66
    & 0.08\\
    bnl2 & 1.1e-04 & 8.20-02 & 0.77 & 0.08 & pilot87 & 1.4e-05 & 6.69-03 &
    1.55 & 0.14\\
    boeing1 & 3.2e-04 & 5.32-02 & 0.64 & 0.08 & pilotnov & 0.0e+00 & 0.00+00 &
    0.00 & 0.08\\
    boeing2 & 6.1e-04 & 1.31-02 & 0.61 & 0.08 & qap8 & 0.0e+00 & 0.00+00 &
    0.00 & 0.07\\
    bore3d & 3.3e-05 & 4.98-03 & 0.03 & 0.08 & qap12 & 0.0e+00 & 0.00+00 &
    0.00 & 0.10\\
    brandy & 5.7e-05 & 1.07-03 & 0.00 & 0.08 & qap15 & 0.0e+00 & 0.00+00 &
    0.00 & 0.10\\
    capri & 0.0e+00 & 1.60-03 & 0.00 & 0.07 & recipe & 0.0e+00 & 1.23-03 &
    0.00 & 0.07\\
    cre-a & 0.0e+00 & 0.00+00 & 0.00 & 0.07 & sc50a & 0.0e+00 & 4.92-03 & 0.02
    & 0.08\\
    cre-b & 0.0e+00 & 0.00+00 & 0.00 & 0.11 & sc50b & 0.0e+00 & 4.69-03 & 0.01
    & 0.07\\
    cre-c & 0.0e+00 & 0.00+00 & 0.00 & 0.07 & sc105 & 0.0e+00 & 4.99-03 & 0.04
    & 0.07\\
    cre-d & 0.0e+00 & 0.00+00 & 0.00 & 0.11 & sc205 & 0.0e+00 & 4.99-03 & 0.34
    & 0.07\\
    cycle & 6.0e-03 & 2.73-01 & 0.80 & 0.08 & scagr7 & 2.1e-05 & 1.43-05 &
    0.00 & 0.08\\
    czprob & 7.1e-05 & 4.99-03 & 0.11 & 0.10 & scagr25 & 1.3e-05 & 1.91-05 &
    0.00 & 0.09\\
    d2q06c & 1.3e-04 & 1.37-02 & 1.20 & 0.14 & scfxm1 & 9.2e-05 & 2.79-03 &
    0.01 & 0.08\\
    d6cube & 0.0e+00 & 0.00+00 & 0.00 & 0.11 & scfxm2 & 9.8e-05 & 4.65-03 &
    0.01 & 0.08\\
    degen2 & 6.3e-04 & 9.83-03 & 0.65 & 0.09 & scfxm3 & 8.0e-05 & 4.93-03 &
    0.01 & 0.09\\
    degen3 & 3.9e-04 & 3.62-03 & 1.00 & 0.09 & scorpion & 0.0e+00 & 0.00+00 &
    0.00 & 0.07\\
    dfl001 & 0.0e+00 & 2.83-03 & 0.09 & 0.11 & scrs8 & 1.6e-03 & 1.88-01 &
    0.65 & 0.08\\
    e226 & 2.5e-02 & 1.00-02 & 0.63 & 0.09 & scsd1 & 0.0e+00 & 0.00+00 & 0.00
    & 0.07\\
    etamacro & 3.9e-04 & 5.43-03 & 0.63 & 0.08 & scsd6 & 0.0e+00 & 0.00+00 &
    0.00 & 0.07\\
    fffff800 & 2.9e-06 & 4.46-03 & 0.54 & 0.09 & scsd8 & 0.0e+00 & 0.00+00 &
    0.00 & 0.07\\
    finnis & 0.0e+00 & 0.00+00 & 0.00 & 0.07 & sctap1 & 0.0e+00 & 0.00+00 &
    0.00 & 0.07\\
    fit1d & 8.3e-01 & 8.90-01 & 0.76 & 0.08 & sctap2 & 0.0e+00 & 0.00+00 &
    0.00 & 0.08\\
    fit1p & 0.0e+00 & 0.00+00 & 0.00 & 0.07 & sctap3 & 0.0e+00 & 0.00+00 &
    0.00 & 0.08\\
    fit2d & 4.2e-02 & 6.04-01 & 2.65 & 0.20 & seba & 0.0e+00 & 0.00+00 & 0.00
    & 0.07\\
    fit2p & 0.0e+00 & 0.00+00 & 0.00 & 0.10 & share1b & 0.0e+00 & 3.22-02 &
    0.61 & 0.08\\
    forplan & 8.6e-05 & 1.32-03 & 0.01 & 0.08 & share2b & 9.5e-04 & 3.47-01 &
    0.60 & 0.08\\
    ganges & 0.0e+00 & 6.45-04 & 0.00 & 0.07 & shell & 0.0e+00 & 0.00+00 &
    0.00 & 0.07\\
    gfrd-pnc & 0.0e+00 & 2.35-16 & 0.00 & 0.07 & ship04l & 5.8e-04 & 5.57-02 &
    0.69 & 0.09\\
    greenbea & 1.2e-05 & 6.78-05 & 0.00 & 0.08 & ship04s & 2.5e-04 & 5.57-03 &
    0.66 & 0.09\\
    greenbeb & 4.4e-05 & 8.66-05 & 0.00 & 0.09 & ship08l & 0.0e+00 & 0.00+00 &
    0.00 & 0.07\\
    grow7 & 0.0e+00 & 0.00+00 & 0.00 & 0.07 & ship08s & 2.4e-04 & 3.72-03 &
    0.68 & 0.08\\
    grow15 & 0.0e+00 & 0.00+00 & 0.00 & 0.07 & ship12l & 0.0e+00 & 0.00+00 &
    0.00 & 0.08\\
    grow22 & 0.0e+00 & 0.00+00 & 0.00 & 0.07 & ship12s & 1.0e-04 & 3.76-03 &
    0.51 & 0.09\\
    israel & 2.0e-05 & 4.50-03 & 0.04 & 0.08 & sierra & 0.0e+00 & 1.87-01 &
    0.72 & 0.09\\
    kb2 & 0.0e+00 & 5.69-02 & 0.63 & 0.08 & stair & 0.0e+00 & 4.31-03 & 0.00 &
    0.08\\
    ken-07 & 2.4e-04 & 4.72-04 & 0.76 & 0.08 & standata & 0.0e+00 & 0.00+00 &
    0.00 & 0.08\\
    ken-11 & 7.8e-05 & 1.95-04 & 0.00 & 0.09 & standgub & 0.0e+00 & 0.00+00 &
    0.00 & 0.07\\
    ken-13 & 1.9e-05 & 6.17-05 & 0.00 & 0.14 & stan. & 0.0e+00 & 0.00+00 &
    0.00 & 0.07\\
    ken-18 & 1.6e-05 & 4.42-05 & 0.00 & 0.18 & stocfor1 & 9.3e-03 & 3.47-03 &
    0.61 & 0.07\\
    lotfi & 0.0e+00 & 3.45-02 & 0.62 & 0.08 & stocfor2 & 9.4e-04 & 5.08-04 &
    0.72 & 0.08\\
    maros-r7 & 0.0e+00 & 0.00+00 & 0.00 & 0.11 & stocfor3 & 2.2e-04 & 1.79-04
    & 1.66 & 0.13\\
    maros & 6.3e-04 & 2.48-03 & 0.67 & 0.08 & truss & 0.0e+00 & 0.00+00 & 0.00
    & 0.07\\
    modszk1 & 0.0e+00 & 0.00+00 & 0.00 & 0.07 & tuff & 0.0e+00 & 0.00+00 &
    0.00 & 0.08\\
    nesm & 0.0e+00 & 0.00+00 & 0.00 & 0.08 & vtp.base & 0.0e+00 & 0.00+00 &
    0.00 & 0.07\\
    osa-07 & 0.0e+00 & 0.00+00 & 0.00 & 0.11 & wood1p & 0.0e+00 & 0.00+00 &
    0.00 & 0.10\\
    osa-14 & 0.0e+00 & 0.00+00 & 0.00 & 0.11 & woodw & 0.0e+00 & 0.00+00 &
    0.00 & 0.07 \\
\bottomrule
\end{tabular}
}
\end{table}


\begin{rem}
  Since our modification of \texttt{Netlib} instances destroys the original bound and coefficients,
  our experiment on \texttt{Netlib} is presented only for reference.
\end{rem}

\subsection{Sifting and Column Generation} \label{app:sifting}

This section discusses online LP in relation to sifting and column generation. We first explain the basic idea and history of the sifting algorithm. Then we identify some common challenges and solutions for sifting solvers. Finally, we show that our methods fit well with sifting and how we can make use of outputs of the online algorithms to speed up sifting.

\subsubsection{Sifting: Column Generation for Linear Programming}

Assume we are solving a standard-form LP
\begin{eqnarray*}
  \max_{\x} & \langle \tmc, \x \rangle & \\
  \text{subject to} & \A \x = \tmb & \\
  & \x \geq \0 & 
\end{eqnarray*}
where $\A \in \mathbb{R}^{m \times n}$ and suppose $n \gg m$. By the theory of
LP, we know that given an optimal basis $\mathcal{B}^{\ast}$ and the
corresponding solution $\x^{\ast}$,
$ | \{ i|x_i^{\ast} > 0 \} | \leq m \ll
   n. $
This property tells us one important fact that {\tmem{most columns are
redundant}} for achieving the optimal solution. One illustration of this property is that if we
somehow manage to find a subset $\mathcal{B}^{\ast} \subseteq \mathcal{W}
\subseteq [n]$ and solve a smaller LP
\begin{eqnarray*}
  \max_{\x_{\mathcal{W}}} & \langle \tmc_{\mathcal{W}}, \x_{\mathcal{W}}
  \rangle & \\
  \text{subject to} & \A_{\mathcal{W}} \x_{\mathcal{W}} = \tmb & \\
  & \x_{\mathcal{W}} \geq \0, & 
\end{eqnarray*}
then $\langle \tmc_{\mathcal{W}}, \x_{\mathcal{W}}^{\ast} \rangle =
\langle \tmc_{\mathcal{W}}, \x_{\mathcal{W}}^{\ast} \rangle +
\langle \tmc_{\bar{\mathcal{W}}}, \x_{\bar{\mathcal{W}}}^{\ast}
\rangle = \langle \tmc, \x^{\ast} \rangle$ and we solve the
original LP at lower cost by concatenating $\x^{\ast} = (
\x_{\mathcal{W}}^{\ast}, \x_{\bar{\mathcal{W}}}^{\ast} = \0 )$. The very
intuition behind sifting is to iteratively update $\mathcal{W}$ till
$\mathcal{B}^{\ast} \subseteq \mathcal{W}$, while in hope of $| \mathcal{W} | \ll
[n]$. \condsep

Although we solve for the optimal primal solution $\x^{\ast}$, one of the most
important components of sifting instead lies in the dual solution $\y^{\ast}$,
as dual solutions tell us how to update $\mathcal{W}$ if $\mathcal{B}^{\ast}
\nsubseteq \mathcal{W}$. Given optimal $( \x_{\mathcal{W}}^{\ast},
\y_{\mathcal{W}}^{\ast} )$ to the aforementioned LP but $\mathcal{W}$
does not contain any optimal bases, we know that $\tmc - \A^{\top}
\y_{\mathcal{W}}^{\ast} \nleq \0$, and otherwise $\x^{\ast} = (
\x_{\mathcal{W}}^{\ast} ; \x_{\bar{\mathcal{W}}}^{\ast} = \0 )$ is
certificated as optimal. In other words, we know some dual infeasibility hides
in $\bar{\mathcal{W}}$:
\[ c_j - \langle \tma_j, \y_{\mathcal{W}}^{\ast} \rangle > 0,
   \text{\quad for some } j \in \mathcal{I} \subseteq \bar{\mathcal{W}} . \]
and this information guides us to eliminate such dual infeasibility by
updating $\mathcal{W} \leftarrow \mathcal{W} \cup \mathcal{I}$. Now that
$(\x_{\mathcal{W}}^{\ast} ; \x_{\bar{\mathcal{W}}}^{\ast} )$ is feasible
for the updated LP problem, sifting subproblems can be efficiently warm-started from
the previous iterations. \condsep

Now we are ready to formalize the LP sifting procedure. Given $\mathcal{W}
\subseteq [n]$, the smaller LP is called a \textit{working problem} and the
procedure finding $\mathcal{I} \subseteq \bar{\mathcal{W}}$ is called
\textit{pricing}. Sifting iteratively updates the working problem by
pricing out dual infeasible columns, till $\mathcal{I}= \varnothing$.

\paragraph{Development of LP Sifting}
The idea of sifting was initially proposed in {\cite{forrest1989mathematical}}
and formalized in a case study {\cite{bixby1992very}} to solve the LP
relaxation of a huge airline crew scheduling problem. One advantage of sifting
is the freedom to choose solvers for the working problems, such as simplex and
the interior point method. Since then, sifting has been adopted as a common
framework for huge LPs and applied in different fields
{\cite{kasirzadeh2017airline,moharrami2015elastoplastic,yanover2006linear}}.
State-of-the-art mathematical programming softwares
{\cite{berthold2018parallelization,pedroso2011optimization,ge2022cardinal,manual1987ibm}}
these days implement their own sifting solvers and trigger it when the ratio
$n / m$ is large. In a word, sifting has evolved into a mature engineering
technique for huge LPs.

\paragraph{Connection with Column Generation}
It's not hard to see sifting bears great resemblance to column generation
{\cite{lubbecke2005selected}} from early LP/MIP literature. To some extent
sifting is a special case of column generation and one major difference lies
in the treatment of the pricing problem. In sifting we can enumerate all the
columns to find $\mathcal{I}= \{ j \in \bar{\mathcal{W}} : c_j -
\langle \tma_j, \y_{\mathcal{W}}^{\ast} \rangle > 0 \}$,
while in the traditional setting of column generation (for example, cutting
stock {\cite{ben2005cutting}}) we often resort to heuristics or combinatorial
approaches to identify ``the most infeasible column''
\begin{eqnarray*}
  \max_{j \in \bar{\mathcal{W}}} & c_j - \langle \tma_j,
  \y_{\mathcal{W}}^{\ast} \rangle . & 
\end{eqnarray*}
Due to the deep connection between sifting and column generation, most
techniques developed for column generation can be smoothly applied to sifting.
In the next section, we discuss the difficulties of implementing a sifting solver and some widely known solutions from the column generation literature.

\subsubsection{Difficulties and Solutions}

In this section, we discuss the difficulties when implementing a sifting
solver. There are three major difficulties {\cite{lubbecke2005selected}} well-known as \textbf{1)}.
heading-in \textbf{2)}. tailing-off and \textbf{3)}. dual oscillation. And
we would like to further ascribe these difficulties to a lack of \textit{prior knowledge}.\condsep

\paragraph{{\underline{Lack of Prior Primal Knowledge}}}

We refer to prior primal knowledge as a measure of \textit{likelihood} that
each column participates in the optimal, or simply a feasible basis. We
believe a lack of this knowledge is partially responsible for the heading-in
and tailing-off effect.\condsep

Heading-in effect appears in the initialization of sifting, where we have to
start from some initial working problem $\mathcal{W}$ and move on. However, if
we are given no prior knowledge, how to pick $\mathcal{W}$ becomes a problem:
it's unlikely that arbitrarily initialized $\mathcal{W}$ would produce a
feasible, not to mention an approximately optimal solution to the original
problem. Therefore most sifting implementations resort to big-$M$ method,
where the original problem is augmented by two blocks of artificial variables
associated with big-$M$ penalties.
\begin{eqnarray*}
  \max_{\x_{\mathcal{W}}, \s_l, \s_u} & \left\langle \tmc_{\mathcal{W}}, \x_{\mathcal{W}}
  \right\rangle - M \left\langle \e, \s_l \right\rangle - M \left\langle \e,
  \s_u \right\rangle & \\
  \text{subject to} & \A_{\mathcal{W}} \x_{\mathcal{W}} - \s_l + \s_u = \tmb &
  \\
  & \x_{\mathcal{W}}, \s_l, \s_u \geq \0 & 
\end{eqnarray*}
The augmented problem is equivalent to the original problem if $M$ is
sufficiently large, and it admits a trivial initial feasible solution. However, whenever $\s_l$ or $\s_u$
has an entry in the basis, $M \left\langle \e, \s_l \right\rangle + M
\left\langle \e, \s_u \right\rangle$ would make the objective value from
sifting far from the true approximate objective, and thus the initial sifting
iterations provide little information about the original problem: we get little
information until we kick $\s_l$ and $\s_u$ out of basis. This effect is known
as heading-in and can be addressed by the prior knowledge about an approximate
feasible primal solution.\condsep

Tailing-off refers to the phenomenon where consecutive sifting iterations
bring little progress when sifting converges. In other words, at the end of
sifting we keep pricing out ``useless'' columns that bring no actual
improvement, and the true optimal basic columns stay in $\bar{\mathcal{W}}$.
While several factors may contribute to tailing-off, prior knowledge about
some approximate optimal solution can efficiently alleviate this effect.
Namely given some approximate optimal solution $\hat{\x}$, we could either
incorporate $\hat{\x}$ in the pricing rule, or simply keep all the basic
columns from $\hat{\x}$ in the working problem.\condsep

\paragraph{{\underline{Lack of Prior Dual Knowledge}}}

We refer to dual prior knowledge as an approximate dual optimal solution. A
lack of this knowledge directly results in the notorious dual oscillation
effect in sifting. For a more rigorous definition and analysis of dual
oscillation, we refer the interested readers to
{\cite{lubbecke2005selected,desrosiers2005primer}}, and in a word dual
oscillation refers to the unstable behavior of the dual sequence $\left\{ \y^{\ast}_{\mathcal{W}}
\right\}$, which is also one important reason why
tailing-off happens. There is vast literature attacking the issue of dual
oscillation, and a systematic approach, known as dual stabilization, has been
proposed and successfully applied to many applications. An important aspect of
dual stabilization is to make $\left\{ \y^{\ast}_{\mathcal{W}} \right\}$ go
more smoothly by taking average
\[ \y^{\ast}_{\mathcal{W}} \leftarrow \alpha \y^{\ast}_{\mathcal{W}} + (1 -
   \alpha) \hat{\y}, \alpha \in (0, 1] \]
where $\hat{\y}$ is an anchor point in the dual space obtained either
before or during sifting. Some common choices of $\hat{\y}$ are geometric
centers of the primal polytope $\left\{ \x : \A \x = \tmb, \x \geq \0
\right\}$, such as analytic center {\cite{luo1998analytic}} and Chebyshev
center {\cite{lee2011chebyshev}}. But these centers are generally too costly
to be computed and less practical for really huge-scale problems.\condsep

So far we have discussed several issues sifting faces and their solutions.
Overall we should find some approximate primal/dual optimal solution at low
cost, and our online algorithms has a role to play here.

\subsubsection{Accelerated Sifting via Online Algorithms}

Finally, we are ready to present our accelerated sifting procedure using
online algorithms. Recall that our method \textbf{1)}. outputs an
approximate primal estimate $\hat{\x}$. \textbf{2)}. outputs a dual
approximate solution. \textbf{3)}. runs in $\mathcal{O} \left( \text{nnz} \left( \A
\right) \right)$ time. Therefore our method can provide both primal and dual
estimates at very low cost.
\begin{figure*}[h]
	\centering
	\fbox{\includegraphics[scale=0.55]{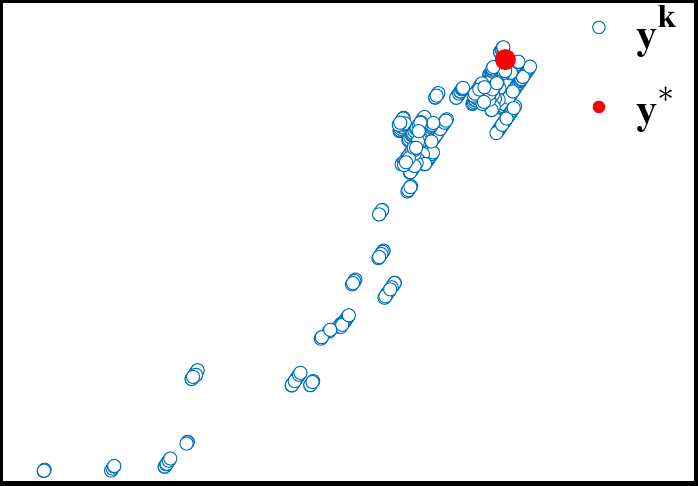}}
		\caption{Convergence of the dual solution $\y^k$ to the optimal $\y^*$ \vspace{10pt}}
\end{figure*}

\begin{rem}
  Generally online algorithm does not address heading-in since the problem
  structure we target admits a trivial primal feasible solution $\x = \0$. But providing an initial 
  good guess of optimal basis can still speed up sifting.
\end{rem}

\subsection{Integer Programming} \label{app:mip}

Finally, we remark that our method can be naturally extended to binary (and
integer) programming. Let $\x_{\text{Bin}}^{\ast}$ denote the optimal solution to the binary
problem and let $\x_{\text{LP}}^{\ast}, \y_{\text{LP}}^{\ast}$ be the optimal
solution to the LP relaxation. Also let $\hat{\x}_{\text{Bin}}$ be some
integer feasible solution and $\hat{\y}_{\text{LP}}$ be some dual feasible
solution, then the following chain of inequalities hold
\[ \langle \tmc, \hat{\x}_{\text{Bin}} \rangle \leq
   \langle \tmc, \x_{\text{Bin}}^{\ast} \rangle \leq \langle
   \tmc, \x_{\text{LP}}^{\ast} \rangle \leq \langle \tmb,
   \y_{\text{LP}}^{\ast} \rangle + \langle \1, [ \tmc -
   \A^{\top} \y_{\text{LP}}^{\ast} ]_+ \rangle \leq \langle
   \tmb, \hat{\y}_{\text{LP}} \rangle + \langle \1, [ \tmc
   - \A^{\top} \hat{\y}_{\text{LP}} ]_+ \rangle \]
and
\[ \langle \tmc, \x_{\text{LP}}^{\ast} \rangle - \langle
   \tmc, \hat{\x}_{\text{Bin}} \rangle \leq \langle \tmb,
   \hat{\y}_{\text{LP}} \rangle + \langle \1, [ \tmc -
   \A^{\top} \hat{\y}_{\text{LP}} ]_+ \rangle - \langle
   \tmc, \hat{\x}_{\text{Bin}} \rangle, \]
which implies $\hat{\y}$ provides a valid dual bound for the binary programming
problem. Combined with $\hat{\x}$ obtained by rounding the solution from online algorithm $\hat{\x}_{\text{LP}}$, we can expect approximately solving a binary
programming problem without resorting to branch and bound. As with general
integer problems, we can split $x_j \in \{ 0, \ldots, U \}$ into $x_j =
\sum_{k = 1}^U x_{j k}$ and re-apply the algorithm for binary problems.


\end{document}